\documentclass[11pt,a4paper,leqno]{article}
\usepackage{amsmath}
\usepackage{amssymb}
\advance\topmargin by -2.2cm
\newcommand{\bbr}{I\!\!R}

\newcommand{\bbn}{I\!\!N}

\newcommand{\bbf}{I\!\!F}
\newcommand{\bbg}{I\!\!\!\!G}

\newcommand{\cala}{{\cal A}}
\newcommand{\calb}{{\cal B}}
\newcommand{\calc}{{\cal C}}

\newcommand{\cale}{{\cal E}}
\newcommand{\calf}{{\cal F}}
\newcommand{\calg}{{\cal G}}

\newcommand{\call}{{\cal L}}

\newcommand{\caln}{{\cal N}}

\newcommand{\calr}{{\cal R}}

\newcommand{\barr}{\begin{array}}
\newcommand{\earr}{\end{array}}
\newcommand{\beqq}{\begin{equation}}
\newcommand{\eeqq}{\end{equation}}
\newcommand{\beao}{\begin{eqnarray*}}
\newcommand{\eeao}{\end{eqnarray*}\noindent}
\newcommand{\beam}{\begin{eqnarray}}
\newcommand{\eeam}{\end{eqnarray}\noindent}

\newcommand{\halmos}{\quad\hfill\mbox{$\Box$}}

\newcommand{\la}{\lambda}

\newcommand{\si}{\sigma}
\newcommand{\al}{\alpha}
\newcommand{\vth}{\vartheta}

\newcommand{\om}{\omega}
\newcommand{\Om}{\Omega}

\newcommand{\vep}{\varepsilon}

\newcommand{\wh}{\widehat}
\newcommand{\wt}{\widetilde}

\newcommand{\ol}{\overline} 
\newcommand{\ov}{\overline} 
\newcommand{\ul}{\underline}

\newcommand{\Lra}{\Longrightarrow}

\newcommand{\lra}{\longrightarrow}

\newcommand{\nto}{n\to\infty}

\begin{document}

{\huge  Estimating discontinuous periodic signals\\ in a time inhomogeneous diffusion}\\

{\bf Reinhard H\"opfner and Yury Kutoyants} 

\vskip0.2cm
Johannes-Gutenberg-Universit\"at, Mainz, and Universit\'e du Maine, Le Mans\\

\vskip1.2cm
{\small 
{\bf Abstract: } 
We consider a diffusion $(\xi_t)_{t\ge 0}$ with some $T$-periodic time dependent input term contained in the drift: under an unknown parameter $\vth\in\Theta$, some discontinuity --~an additional periodic  signal~--  occurs at times $kT{+}\vth$, $k\in\bbn$. Assuming positive Harris recurrence of $(\xi_{kT})_{k\in\bbn _0}$ and exploiting the periodicity structure, we prove limit theorems for certain martingales and functionals of the process $(\xi_t)_{t\ge 0}$.  
They allow to consider the statistical model parametrized by $\vth\in\Theta$ locally in small neighbourhoods of some fixed $\vth$, with radius $\frac{1}{n}$ as $\nto$. We prove convergence of local models to a limit experiment studied by Ibragimov and Khasminskii [IH~81] and discuss the behaviour of estimators under contiguous alternatives. 

likelihood ratio processes, convergence of experiments, contiguity, maximum likelihood estimators, Bayes estimators, local asymptotic minimax bound. 

\vskip0.2cm
{\bf MSC:  \quad 62 F 12 ,  60 J 60  }
}%%ende von small

\vskip1.2cm
We consider a problem of parameter estimation in a Markov process $\left(\xi_t\right)_{t\ge 0}$ whose drift is $T$--periodic in the time variable. A  parameter $\vth$ comes in through a periodic signal with periodicity $T$, and represents a time of discontinuity. Our main assumption on the process is positive Harris recurrence of $\left(\xi_{kT}\right)_{k\in\bbn _0}$ from which --~exploiting the periodicity structure of the semigroup~-- we prove positive Harris recurrence of the chain of $T$--segments $\left(\, (\xi_{kT+s})_{0\le s\le T} \,\right)_{k\in\bbn _0}$. This allows to deduce limit theorems for certain martingales and strong laws of large numbers for certain functionals of the continuous-time process $(\xi_t)_{t\ge 0}$. Based on these, we deal with convergence of local models at $\vth$ --~corresponding to observation of the process $\left(\xi_t\right)_{t\ge 0}$ up to time $nT$, local scale at $\vth$ will be $\frac{1}{n}$ as $\nto$~--  to a limit model whose likelihoods are of type $u \to e^{W_u-\frac12|u|}$ with double-sided Brownian motion $W$. In this limit model, investigated by Ibragimov and Khasminskii [IH 81, sections VII.2--3], parametrization is $\frac12$--H\"older continuous in the sense of Hellinger distance. So we are far from the framework of local asymptotic (mixed) normality (well-studied since LeCam [L 68] and H\'ajek [H 70], see [D 85], [LY 90]), or from $L^2$-differentiable experiments. Ibragimov and Khasminskii considered a simple 'signal in white noise' setting, proved convergence of maximum likelihood and Bayes estimators at $\vth$, calculated the limit variance of the maximum likelihood estimator, and pointed out that in a limit model with likelihood ratios of type  $u \to e^{W_u-\frac12|u|}$ a Bayes estimator is better than the maximum likelihood estimator. 

\vskip0.3cm
Convergence to this limit model has been considered in several settings since then: see [KK~00] in a framework of delay equations,  [K~04, section 3.4] in time homogeneous ergodic diffusions with spatial discontinuity in the drift, [DP 84, section 3] in the context of a change point in iid observations, [D~09] and [CK~09] with two-sided compensated Poisson or compound Poisson processes appearing in log-likelihood ratios.  The approach of Ibragimov and Khasminskii starts from certain assumptions on Hellinger distances and from convergence of likelihood ratios 'uniformly in $\vth$'.  In several aspects, our approach is different. We develop limit theorems in diffusions with $T$--periodic semigroup  which will be our key tool in view of  convergence of likelihood ratios and estimators. Whereas these  allow to check and exploit assumptions on Hellinger distances similiarly to the work quoted above, convergence of likelihood ratios 'uniformly in $\vth$' is not suitable for our framework of inhomogeneity in time, and is systematically avoided.  We are focussing on contiguous alternatives, make extensive use of  'LeCam's Third Lemma' (see [LY 90, pp.\ 22--23]), and exploit asymptotic equivariance of suitable estimator sequences with respect to contiguous alternatives.  Our local asymptotic minimax bound, of the type of the asymptotic minimax bound for risk functions in Strasser ([S~85, Cor.\ 62.6], see also [L~72], [M~82], [V~91]), controls a maximal quadratic risk on shrinking neighbourhoods of $\vth$ with radius proportional to $\frac{1}{n}$, and a Bayes sequence attains this bound.

\vskip0.3cm 
We describe our setting in more detail. The observed diffusion process is inhomogeneous in time 
\beqq\label{process}
d\xi_t  \;=\;  \left[ S(\vth,t) + b(\xi_t) \right]dt   \;+\;  \si(\xi_t)\,dW_t  \;,\quad t\ge 0 \;. 
\eeqq
with some deterministic and $T$-periodic input  $t\to S(\vth,t)$ where, for known periodicity $T$ and for 
known functions $\la(\cdot)\ge 0$ and $\la^*(\cdot)>0$ which are continuous and $T$-periodic,  
\beqq\label{signal}
S(\vth,t) \;=\;  \la(t) \;+\; \la^*(t)\, 1_{(\vth,\vth+a)} (i_T(t))  \;,\; t\ge 0 \;,\quad\mbox{with}\;\;  i_T(t) \;:=\; t \;\mbox{modulo}\; T   
\eeqq
depends on an unknown parameter $\vth$. This means that some additional input $\la^*$ is switched on periodically at times $kT{+}\vth$, $k\in\bbn_0$, and is of known duration $a$. We put $\Theta:=(0,T-a)$. \\
The functions  $b(\cdot)$ and  $\si(\cdot)$ are Lipschitz; hence for all values of the parameter $\vth$,  we have Lipschitz and linear growth conditions for the time-dependent coefficients of the above SDE, and thus existence  and pathwise uniqueness for its solution. 

\vskip0.3cm
We are interested in convergence of local models and convergence of maximum likelihood (MLE) and Bayes (BE) estimators for the unknown parameter $\vth\in\Theta$ when a trajectory of $\xi$ has been observed up to time $nT$. As $\nto$, the right choice of local scale for local models at $\vth$ turns out to be $\frac{1}{n}$, at every point $\vth\in\Theta$. For the limit  of local models at $\vth$, we find  likelihood ratios  
\beqq\label{limitlikelihoodratiogeneral}
\wt L^{ u / 0 } \;:=\;  \exp \left\{\;   \wt W( u J_\vth )  -  \frac{1}{2} \left| u J_\vth \right|  \;\right\}  \;,  \quad u\in\bbr
\eeqq
with double-sided Brownian motion $(\wt W_u)_{u\in\bbr}$ and with scaling constants  $0<J_\vth<\infty$. From Terent'yev [T 68] over Golubev [G 79] and Ibragimov and Khasminskii [IH 81]  to Rubin and Song [RS 95] it has become evident that in experiments of type  (\ref{limitlikelihoodratiogeneral}) with unknown parameter $u\in\bbr$, the variance of the Bayes estimator with respect to quadratic loss is strictly smaller than the variance of the maximum likelihood estimator. 
Note that in the limit experiment (\ref{limitlikelihoodratiogeneral}) --~not a 'quadratic' experiment, not $L^2$-differentiable~--  there is no sufficient statistic,  
no analogue of a central statistic, and no analogue of a convolution theorem ([H~70], [J 82]) which in the classical LAN/LAMN case--~together with a lemma stating that 'arbitrary estimator sequences are in some sense almost equivariant'~--  is the key tool to obtain a local asymptotic minimax theorem (see [LY 90, p.\ 83]) simultaneously for a large class of loss functions. 
In the context here, we have not more than local asymptotic minimax bounds with respect to specified loss functions (see Strasser [S~85, Cor.\ 62.6], or [L~72], [M~82], [V~91]), and no tool for comparison between different bounds associated to different choices of a loss function. 

Using contiguity techniques and exploiting an equivariance property (lemma 5.3 below) of the limit experiment (\ref{limitlikelihoodratiogeneral}), we consider squared loss and prove a bound (theorem 1.8 below) 
\beqq\label{LAMAnnounce}
\lim_{C\uparrow\infty}\; \liminf_{\nto}\; \inf_{\wt \vth_{nT}}\; \sup_{|u|\le C}\; 
E_{\vth+\frac{u}{n}}\left(\; \left[ n\left( \wt \vth_{nT} - (\vth{+}\frac{u}{n}) \right) \right]^2 \;\right) 
\quad\ge\quad E\left( \left[ u^* \right]^2 \right)
\eeqq
where $\,u^*$ is the Bayes estimator associated to quadratic loss in the limit experiment (\ref{limitlikelihoodratiogeneral}), and where $\inf\limits_{\wt \vth_{nT}}$  allows to compare all possible estimators based on observation of $\xi$ up to time $nT$. We prove asymptotic equivariance with respect to contiguous alternatives of a Bayes sequence $\vth^*_{nT}$ with respect to squared loss (proposition 1.7.b) as $\nto$: thus this Bayes sequence attains the bound (\ref{LAMAnnounce}). 

\vskip0.3cm
In view of asymptotic statistical properties,  our model  behaves exactly as the simple 'signal in white noise' setting of [IH 81, section VII.2] which corresponds to the special case $\si(\cdot)\equiv 1$, $b(\cdot)\equiv 1$ in (\ref{process}) above. However, if limit theorems are the key tool to prove statistical properties (convergence of experiments, convergence of estimators, ...), these are radically different in our case. In their  likelihoods, thanks to  $\si(\cdot)\equiv 1$,  [IH~81] can work with very simple Gaussian processes where calculation of means and covariances is enough to determine the limiting behaviour. Similiarly,  in the time homogeneous ergodic diffusion model of [K~04, section 3.4] with one discontinuity in the drift, well known limit theorems for convergence of martingales and of additive functionals for ergodic diffusions are at hand. Our time inhomogeneous $T$--periodic problem (\ref{process}) with non-trivial $\si(\cdot)$ requires a completely new approach.  So an essential  part of the present paper is devoted to proving limit theorems which make statistical theories work in our setting. Also in view of the behaviour of estimators under contiguous alternatives, we have to go beyond what had been done earlier in order to obtain the local asymptotic minimax bound (\ref{LAMAnnounce}).  

\vskip0.3cm
Our interest in periodicity structures in diffusions is linked to the following application. In some membrane potential data sets similiar to those investigated in [H 07] which we wish to interprete as realizations of certain SDE's (out of many references, we mention [LL 87], [T 89], [LS 99], [DL 05], [DL 06]), there is evidence for time-dependent  'input' in the drift which the modelization has to take into account. In analogy to the result of [BH~06, section 3.2] on  large systems of neurons receiving identical time-dependent input,  questions of  periodic input received by a single neuron in an active network deserve to be studied. In particular, a discontinuity (\ref{signal}) with constants $\la$, $\la^*$  can be interpreted as some stimulus switched on/off periodically, and is of biological relevance. 

\vskip0.3cm
This paper is organized as follows. Section 1 states all statistical results concerning the model defined by (\ref{process})+(\ref{signal}): convergence of experiments, convergence of estimators, local asymptotic minimax bound (theorem 1.8).  Section 2 deals with Harris properties of the chain of $T$--segments; it is formulated in a more general setting  and can be read independently. Section 3  contains an exponential inequality adapted to our purposes from Brandt [B~05]). Section 4 works out the limit theorems which we need to consider local models in our problem (\ref{process})+(\ref{signal}). The main results in this probabilistic part of the paper are theorems 2.1+4.1 (strong laws of large numbers for time inhomogeneous diffusions with periodicity structure) and theorem 4.3 with remark 4.4 (convergence of martingale terms which occur in the log-likelihood ratios of  local models at $\vth$). On this basis, section 5 contains the statistical part of work to be done, and collects all proofs for the results stated in section 1.  
\\

%%%%%%%%%%%%%%%%%%%%%%%%%%%%%%%%%%%%%%%%%%%%%
%\newpage
\section*{1~ Outline of statistical results} 
%%%%%%%%%%%%%%%%%%%%%%%%%%%%%%%%%%%%%%%%%%%%%

In order to exploit 'ergodicity properties' of the time inhomogeneous diffusion $\xi = (\xi_t)_{t\ge 0}$ in (\ref{process}) with $T$--periodic time-dependent input (\ref{signal}), our principal assumption will be 
$$
\mbox{the embedded chain $(\xi_{kT})_{k\in\bbn_0}$ is positive recurrent in the sense of Harris} 
\leqno{(H1)}
$$
for any fixed value of the parameter $\vth\in\Theta$. As an example, $(H1)$ always holds if piecewise continuous $T$-periodic input is added to an Ornstein-Uhlenbeck SDE, see 2.3 below. Under $(H1)$, there is a unique invariant probability $\mu^{(\vth)}$ for the chain $(\xi_{kT})_{k\in\bbn_0}$ under $\vth$. We introduce the chain of $T$--segments in the path of $\xi$ 
$$
X = (X_k)_k  \quad\mbox{defined by}\quad  X_k \;:=\; (\xi_{(k-1)T+s})_{0\le s\le T} \;,\; k\ge 1
$$
which --~as a consequence of the $T$--periodicity in the drift of our SDE~-- is time homogeneous. $X$ takes values in the path space $(C_T,\calc_T)$ of continuous functions $[0,T]\to\bbr$. We deduce  from assumption (H1) --~see theorem 2.1 in section 2 below~-- that this $T$--segment chain is positive recurrent in the sense of Harris under $\vth$, with a specified invariant probability $m^{(\vth)}$ on $(C_T,\calc_T)$. Limit theorems for functionals of the process $\xi$ which we need for our analysis of the statistical model  (see theorems 4.1 and 4.3 in section 4 below) are then obtained through strong laws of large numbers in the Harris chain $X = (X_k)_{k\in\bbn_0}$. \\

Let $Q^\vth$ denote the law of the process  $\xi = (\xi_t)_{t\ge 0}$ of (\ref{process})+(\ref{signal}) under $\vth\in\Theta$, a law on the canonical path space $(C,\calc)$ of continuous functions $[0,\infty)\to\bbr$ equipped with its canonical filtration $\bbg$.  Our second major assumption  
$$
\mbox{$\si(\cdot)$ is bounded away from $0$ and $\infty$ on $\bbr$}  
\leqno{(H2)}
$$
guarantees that  for any pair of different values $\zeta '\neq\zeta$ in $\Theta$, the laws $Q^\zeta$ and $Q^{\zeta '}$ are locally equivalent with respect to $\bbg$, and the likelihood ratio process of $Q^{\zeta '}$ to  $Q^\zeta$ is 
\beao
L^{\zeta ' / \zeta}_t   
&=&   \exp \left\{ \int_0^t \frac{S(\zeta ',s)-S(\zeta,s)}{\si^2(\eta_s)}\, dM^{\zeta}_s \;-\; \frac12  \int_0^t   \frac{(S(\zeta ',s)-S(\zeta,s))^2}{\si^2(\eta_s)}\,  ds \right\} \;,\quad t\ge 0 
\eeao
where we write $\eta = (\eta_t)_{t\ge 0}$ for the canonical process on  $(C,\calc)$, and $M^{\zeta}$ for the $(Q^\zeta,\bbg)$--martingale part of  $\eta$ (see [JS 87], [LS 81], [Ku 04]). Let $B$ denote a version of  $\int_0^\cdot \frac{1}{\si(\eta_s)}dM^{\zeta}_s$ under the reference point $\zeta\in\Theta$:  $\,B$ is a $\bbg$--Brownian motion under $Q^\zeta$. 
Since $S(\vth,\cdot)$ is --~up to the continuous functions $\la(\cdot)$, $\la^*(\cdot)$, and up to $T$--periodic continuation~-- the indicator function $1_{(\vth,\vth+a)}$, the likelihood ratio $L^{\zeta ' / \zeta}_t$ takes for $\zeta '$ sufficiently close to $\zeta$ the simple form 
\beao
&&\exp \left\{\;  
 \left[ -\;\int_0^t \frac{\la^*(s)}{\si(\eta_s)} \, 1_{(\zeta,\zeta ')}(i_T(s))\, dB_s   \;-\; \frac12  \int_0^t   \left(\frac{\la^*(s))}{\si(\eta_s)}\right)^2\, 1_{(\zeta,\zeta ')}(i_T(s))\, ds \right] \right. \\
&&+\; \left. \left[ \int_0^t \frac{\la^*(s)}{\si(\eta_s)} \, 1_{(\zeta+a,\zeta '+a)}(i_T(s))\, dB_s   \;-\; \frac12  \int_0^t   \left(\frac{\la^*(s))}{\si(\eta_s)}\right)^2\, 1_{(\zeta+a,\zeta '+a)}(i_T(s))\, ds \right] \; \right\}
\eeao
in case $\zeta<\zeta '<\zeta+a$; the same holds with $\zeta$, $\zeta'$ in the intervals and sign $\pm$ in front of the stochastic integrals interchanged if $\zeta'<\zeta<\zeta'+a$. Our first main result --~very easy to see in the 'trivial' case where $\si(\cdot)$ is constant (cf.\ [IH 81], lemma 2.4 on p.\ 334), essentially more difficult for Lipschitz functions $\si(\cdot)$ satisfying $(H1)$ where the argument has to go back to the chain of $T$--segments in the path of $\xi$ (cf.\ theorems 4.1, 4.3 and remark 4.4 in section 4 below)~-- is the following.   \\

{\bf 1.1 Theorem: } Under Lipschitz and linear growth conditions on $b(\cdot)$ and $\si(\cdot)$, under $(H1)$ and $(H2)$, the following holds for every $\vth\in\Theta$: 

a)~we have convergence under $Q^\vth$ as $\nto$ of 
$$
\left(\; L_{nT}^{(\vth + \frac{u}{n}) / \vth } \;\right) _{u \in \Theta_{\vth,n}}  
\quad,\quad  \Theta_{\vth,n} := \{ u\in\bbr : \vth + \frac{u}{n} \in\Theta\}
$$
in the sense of finite dimensional distributions to 
\beqq\label{limitlikelihoodratio}
\wt L = \left(\wt L^{ u / 0 }\right)_ {u\in\bbr}  
\quad,\quad   \wt L^{ u / 0 } \;:=\;  \exp \left\{\;   \wt W(u J_\vth)  -  \frac{1}{2} | u J_\vth |  \;\right\} 
\eeqq
where $(\wt W_u)_{u\in\bbr}$ is two-sided standard Brownian motion, and  $J_\vth$ the scaling constant 
$$
J_\vth  \;:=\;  \left\{  (\la^*(\vth))^2\, (\mu^{(\vth)} P^{(\vth)}_{0,\vth})  + ( \la^*(\vth{+}a))^2\,  (\mu^{(\vth)} P^{(\vth)}_{0,\vth+a}) \right\} (\frac{1}{\si^2})  \;. 
$$
Here  $(P^{(\vth)}_{s,t})_{0\le s<t<\infty}$ denotes the semigroup of the process $(\xi_t)_{t\ge 0}$ under $\vth$, 
$\mu^{(\vth)}$ the invariant measure for $(\xi_{kT})_{k\in\bbn_0}$  according to $(H1)$, and we write for short $\wt\mu(f)$ for $\int f d\wt\mu$. 

b)~ Let $\wt W$ in (\ref{limitlikelihoodratio}) be defined  on some  $(\wt \Om,\wt \cala,\wt P_0)$. Then  $(\wt \Om,\wt \cala)$ carries a limit experiment  
\beqq\label{limitexperiment}
\wt \cale  \;=\; \left\{ \wt P_u : u\in\bbr \right\}  \;\;\mbox{defined by}\;\;  d\wt P_u := \wt L^{ u / 0 } d\wt P_0  \;\;\mbox{on}\;\;  (\wt \Om,\wt \cala)  
\eeqq 
such that we have convergence of experiments: local experiments at $\vth$
 $$
\cale_n^{(\vth)}  \;:=\;   \left\{ Q^{\vth+\frac{u}{n}} \mid \calg_{nT} \,:\; u \in \Theta_{\vth,n}  \right\} \;,\quad n\ge 1
$$
converge  as $\nto$  to the limit experiment $\cale$  in the sense of Strasser ([S 85], p.\ 302). \\

We will prove this theorem, based on the results of sections 2 to 4, in  5.2 below.\\

{\bf 1.2 Remark: } We give an interpretation for the type of limit experiment in (\ref{limitexperiment}), putting $J_\vth=1$ for short.  Recall that two-sided standard Brownian motion $(\wt W_u)_{u\in\bbr}$ on  $(\wt \Om,\wt \cala,\wt P_0)$ means that  two independent  standard Brownian motions $(\wt W ^+_v)_{v\ge 0}$ and $(\wt W ^-_v)_{v\ge 0}$ exist on  $(\wt \Om,\wt \cala,\wt P_0)$ such that $\wt W_u $ is given by $\wt W ^+_u$ if $u\ge 0$, and by $\wt W ^-_{|u|}$ if $u\le 0$. Define   
\beqq\label{interpretationPu}
\wt P_u \;:=\; \left\{ \begin{array}{ll}
\call \left(   \left( \wt W ^+_v + v\wedge u \;,\;   \wt W ^-_v  \right)_{v\ge 0}   \mid \wt P _0 \right) & \mbox{in case}\;\; u\ge 0  \;, \\
\call \left(   \left( \wt W ^+_v  \;,\;  \wt W ^-_v + v\wedge |u|  \right)_{v\ge 0}   \mid \wt P _0 \right) & \mbox{in case}\;\; u\le 0  \;. 
\end{array} \right.
\eeqq
Here constant drift $1$ is added to $\wt W ^+$ in case $u>0$ and to $\wt W ^-$ in case $u<0$, and is switched off at time $|u|$. Consider the likelihood ratio process of $\wt P_u$ to $\wt P_0$ relative to the filtration generated by the bivariate canonical  process on $C([0,\infty),\bbr^2)$ (see proof of lemma 5.3 below, or use  [JS 87], [LS 81]).  If we are allowed to observe over the {\em infinite} time interval $[0,\infty)$, we end up  with the likelihood ratios  
$$
\wt L^{ u / 0 } \;=\;  \exp \left\{\;   \wt W(u)  -  \frac{1}{2} | u |  \;\right\}  
$$
given in (\ref{limitlikelihoodratiogeneral}).  Up to the scaling factor $J_\vth$, this is  the situation of theorem 1.1.\\

%A minor extension of theorem 1.1 --~exactly analogous to [IH 81, p.\ 335]~-- is\\
%
%{\bf 1.3 Remark: } If we replace  $S(\vth,t) \;=\;  \la(t) + \left(\la^* 1_{(\vth,\vth+a)} \right)(i_T(t)) $ defined in (\ref{signal}) by 
%$$
%S(\vth,t) \;=\; \left( \la+ \sum_{i=0}^\ell \la^*_{i+1} 1_{(\vth+a_i,\vth+a_{i+1})} \right)(i_T(t)) 
%$$
%where $0=a_0<a_1<\ldots<a_{\ell+1}=a$ are fixed and known, where $\la^*_i$ are known continuous and strictly positive functions on $[0,T]$, for $1\le i\le \ell$, and  $\la^*_0\equiv 0\equiv \la^*_{l+2}$, then  theorem 1.1 remains valid  with new definition 
%$$
%J_\vth  \;:=\;   \left\{  \sum_{i=0}^{\ell+1} [(\la^*_{i+1}-\la^*_i)(\vth+a_i)]^2\, (\mu^{(\vth)} P^{(\vth)}_{0,\vth+a_i})  \right\} (\frac{1}{\si^2}) \;. 
%$$
%of the scaling constants. This is a simple variant of the proof of theorem 1.1. \halmos\\

{\bf 1.3 Remark: } Instead of fixed duration $a$ of the additional signal $\la^*$ as in (\ref{signal}), we might work with two-dimensional parameter $(\vth_1,\vth_2)$ such that $0<\vth_1<\vth_2<T$, and with 
$$
S(\vth,t) \;=\;  \la(t) \;+\; \la^*(t)\, 1_{(\vth_2,\vth_2)} (i_T(t))  \;\;,\;\; t\ge 0 \;.  
$$
All results which we give here generalize to this two-dimensional problem: the limit experiment will be a product of two models as in 1.2, with scaling constant $(\la^*(\vth_i))^2\, (\mu^{(\vth_1,\vth_2)} P^{(\vth_1,\vth_2)}_{0,\vth_i}) (\frac{1}{\si^2})$ in factor $i$. \\

Based on theorem 1.1,  we proceed using theorems 19--21 in appendix A I.4 of Ibragimov and Khasminskii [IH 81] (up to the uniform in $\vth$ part which can be omitted there both in the assertion and in the assumption of the theorem) and techniques developped in [IH 81,  chapter~I, theorems 5.1+5.2] to show\\

{\bf 1.4 Theorem: } Under Lipschitz and linear growth conditions on $b(\cdot)$ and $\si(\cdot)$, under $(H1)$ and $(H2)$, we have the following for every $\vth\in\Theta$:  

a)~ For $K<\infty$ arbitrarily large,  the likelihood ratios in local models $\cale_n^{(\vth)}$ at $\vth$ 
$$
\left(\; L_{nT}^{(\vth + \frac{u}{n}) / \vth } \;\right) _{u \in [{-}K,K]}  \quad\mbox{under $Q^\vth$}
$$
converge weakly in $C([{-}K,K])$ %(the space of continuous functions $[{-}K,K]\to\bbr$ with the metric of uniform convergence) 
as $\nto$ to the likelihoods 
 $$
\left(\;  \wt L^{ u / 0 } \;\right) _{u \in [{-}K,K]}    \;\mbox{under $\wt P_0$} 
%\;, \quad   \wt L^{ u / 0 }  \;=\;  e^{\;   \wt W(u J_\vth)  -  \frac{1}{2} | u J_\vth | }  
$$
in the limit model $\wt\cale$ of (\ref{limitlikelihoodratio})+(\ref{limitexperiment}). 

b)~  For arbitrary  $p\in\bbn_0$ and $K_0>0$, there are constants $b_1(p,K_0)$ and $b_2$ such that  
\beqq\label{IHconditon-6}
Q^\vth\left( \sup_{ |u| > K }\, |u|^p\, L_{nT}^{(\vth + \frac{u}{n}) / \vth } \;\ge 1\; \right)  \;\;\le\;\; b_1(p,K_0)\, e^{\,-\; b_2\, K} 
\quad,\quad K_0\le K<\infty 
\eeqq
holds for all $\vth\in\Theta$ and all $n\ge 1$,  together with 
\beqq\label{IHconditon-16}
\wt P_0 \left( \sup_{ |u| > K }\, |u|^p\, \wt L^{ u / 0 } \;\ge 1\; \right)  \;\;\le\;\; b_1(p,K_0)\, e^{\,-\; b_2\, K} 
\quad,\quad K_0\le K<\infty \;. 
\eeqq
Here $b_1(p,K_0)$ does not depend on $\vth$ or $n$, and $b_2$ does not depend on $\vth$, $n$, $p$, $K_0$. 

c)~  For arbitrary  $p\in\bbn_0$ and $K_0>0$, there are constants $\wt b_1(p,K_0)$ and $\wt b_2$ such that   
\beqq\label{IHconditon-86}
E_\vth\left(\,  \int_{ \Theta_{\vth,n} \cap \{ |u| > K \} }  |u|^p\, \frac{ L_{nT}^{(\vth + \frac{u}{n}) / \vth } }{ \int_{\Theta_{\vth,n}} L_{nT}^{(\vth + \frac{u'}{n}) / \vth } du' }\; du 
\;\right)  \;\;\le\;\; \wt b_1(p,K_0)\, e^{\,-\; \wt b_2\, K} 
\quad,\quad K_0\le K<\infty 
\eeqq
holds for all $\vth\in\Theta$ and all $n\ge 1$,  together with 
\beqq\label{IHconditon-96}
E_{\wt P_0} \left(\,  \int_{ \{ |u| > K \} }  |u|^p\,  \frac{ \wt L^{ u / 0 } }{ \int_{\bbr}  \wt L^{ u' / 0 } du' }\; du  \;\right)  \;\;\le\;\; \wt b_1(p,K_0)\, e^{\,-\; \wt b_2\, K} 
\quad,\quad K_0\le K<\infty \;. 
\eeqq
Again $\wt b_1(p,K_0)$ does not depend on $\vth$ or $n$, and $\wt b_2$ does not depend on $\vth$, $n$, $p$, $K_0$.

\vskip0.8cm 
The proof of theorem 1.4 will be given in 5.5 below. As consequences of theorem 1.4, we will obtain convergence of maximum likelihood (MLE) and Bayes estimators (BE) for the unknown parameter when a trajectory of $\xi$ is observed up to time $nT$, $\nto$. The MLE sequence is 
\beqq\label {MLEstage-n}
\wh\vth_{nT} \;=\;\mathop{\rm argmax}\limits_{\zeta\in\ov{\Theta}}\, L_{nT}^{\zeta / \zeta_0 } 
\;:=\; \min\left\{ \zeta\in\ov{\Theta} \;:\; L_{nT}^{\zeta / \zeta_0 } = \max_{\zeta'\in\ov{\Theta}}L_{nT}^{\zeta' / \zeta_0 } \right\} 
\quad,\quad n\in\bbn  
\eeqq
with $\zeta_0\in\Theta$ some fixed point, and $\ov{\Theta}$ the closure of $\Theta$. Presence of  'min' in (\ref{MLEstage-n}) guarantees for a measurable selection whenever the argmax is not unique. The MLE in the limit experiment  
\beqq\label {MLElimitexperiment}
\wh u   \;=\;  \mathop{\rm argmax}\limits_{u\in\bbr}\,  \wt L^{ u / 0 } 
\eeqq
is finite--valued and uniquely determined almost surely, 
%since (\ref{IHconditon-16}) of theorem 1.4 controls the decrease of $u \to \wt L^{ u / 0 } $ as $u$ tends to ${+}\infty$ or ${-}\infty$, and 
by [IH 81, chapter VII, lemma 2.5]. The control (\ref{IHconditon-16}) guarantees in particular that a BE 'with uniform prior on the entire real line' 
\beqq\label{BElimitexperiment}
u^*  \;:=\;   \frac{ \int_{-\infty}^\infty u\,  \wt L^{ u / 0 }\, du }{ \int_{-\infty}^\infty  \wt L^{ u / 0 }\, du } 
\eeqq
(sometimes called Pitman estimator) is well defined in the limit experiment. Correspondingly, we consider the BE sequence with uniform prior on $\Theta = (0,T{-}a)$    
\beqq\label{BE}
\vth^*_{nT}  \;:=\;   \frac{ \int_{\Theta} \,\zeta\;\, L_{nT}^{\zeta / \zeta_0 }\; d\zeta }{ \int_{\Theta}  \,L_{nT}^{\zeta / \zeta_0 }\; d\zeta } 
\quad,\quad n\in\bbn  
\eeqq
for the unknown parameter $\vth\in\Theta$, based on observation of  $\xi$ up to time $nT$ (in (\ref{BE}), we might as well use  smooth and strictly positive prior densities on $\Theta$,  but  this generalization turns out to be without interest in view of theorems 1.7 and 1.8 below). \\

{\bf 1.5 Theorem: } Under Lipschitz and linear growth conditions on $b(\cdot)$ and $\si(\cdot)$, under $(H1)$ and $(H2)$, we have the following properties of the MLE and the BE sequence, for every $\vth\in\Theta$: 

a)~ weak convergence  as $\nto$: 
$$
\call \left(\, n\, ( \wh\vth_{nT} - \vth) \mid Q^\vth \,\right)  \;\;\lra \;\;   \call \left(\, \wh u \mid \wt P _0 \,\right) \;, 
$$
$$
\call \left(\, n\, ( \vth^*_{nT} - \vth) \mid Q^\vth \,\right)  \;\;\lra \;\;   \call \left(\, u^* \mid \wt P _0 \,\right)    \;;  
$$

b)~ finite moments of arbitrary order $p\in\bbn$,  and convergence of moments as $\nto$: 
$$
E_\vth \left(\, \left[ n\, ( \wh\vth_{nT} - \vth) \right]^p  \,\right)  \;\;\lra \;\;   E_{\wt P_0} \left(\, | \wh u |^p \,\right) \;, 
$$  
$$
E_\vth \left(\, \left[ n\, ( \vth^*_{nT} - \vth) \right]^p  \,\right)  \;\;\lra \;\;   E_{\wt P_0} \left(\, | u^* |^p \,\right)  \;. 
$$
\vskip0.8cm

Theorem 1.5 will be proved in 5.6 and 5.7 below. Note that part a) of the theorem would be enough for convergence of risks with respect to loss functions which are continuous, subconvex and bounded, i.e.\ 'bowl-shaped' in the sense of LeCam and Yang [LY 90, p.\ 82], whereas part~b) is useful e.g.\  for quadratic loss  to  be considered below. \\

{\bf 1.6 Remark: } For the limit experiment in the setting of remark 1.2 where $J_\vth$ equals $1$, the following is known. Ibragimov and Khasminskii [IH 81, p.\ 342] calculated the variance of the MLE, whereas Rubin and Song [RS 95] calculated the variance of the BE: the ratio of MLE to BE variance is  $26$ to $16\cdot \zeta(3) \approx 19.3$ where $\zeta(s) = \sum_{n=1}^\infty \frac{1}{n^s}$ is Riemann's zeta function. \halmos\\

In our setting  with scaling factor  $J_\vth$ in the limit experiment (\ref{limitlikelihoodratio})+(\ref{limitexperiment}), the results quoted in remark 1.6 yield the following:  the limit law for rescaled MLE errors has finite variance  
\beqq\label{MLElimitvariance}
J_\vth^{-2}\cdot 26
\eeqq
at $\vth$, different from the limit law for rescaled BE errors whose variance is 
\beqq\label{BElimitvariance}
J_\vth^{-2}\cdot 16\cdot \zeta(3)  \;\approx\;  J_\vth^{-2}\cdot 19.3 \;. 
\eeqq
Moreover, by definition of a Bayes estimator with respect to quadratic loss --~this is the $L^2$-projection property of conditional expectations~-- we have at every finite stage $n$ optimality of BE in the sense of integrated risk under quadratic loss functions: comparing all possible estimators $\wt\vth_{nT}$ based on observation of the trajectory of $\xi$ up to time $nT$, a minimum 
\beqq\label{BEoptimality-1}
\min\left\{\,  \int_\Theta E_\vth\left( \left[  n( \wt\vth_{nT} - \vth) \right]^2 \right) d\vth \;\mid\; \mbox{$\wt\vth_{nT}$ is $\calg^\xi_{nT}$--measurable} \,\right\}  
\eeqq
exists and is realized by the integrated risk of the BE (\ref{BE}) with respect to quadratic loss 
\beqq\label{BEoptimality-2}
\int_\Theta E_\vth\left( \left[  n( \vth^*_{nT} - \vth) \right]^2 \right) d\vth   \;, 
\eeqq
for every fixed $n$. This is an elementary and  pre-asymptotic argument, averaging over the whole parameter space  $\Theta$, and does not contain much information on the behaviour of our estimator in small neighbourhoods of a parameter value  $\vth$. The following two results deal with shrinking neighbourhoods of fixed points in the parameter space, in the sense of contiguous alternatives, and fill this gap.   \\

{\bf 1.7 Theorem: } Under Lipschitz and linear growth conditions on $b(\cdot)$ and $\si(\cdot)$, under $(H1)$ and $(H2)$, the following holds  for every $\vth\in\Theta$. 

a)~ The BE $u^*$  in the limit experiment $\wt\cale\,$ with 'uniform prior over $\bbr$' (\ref{BElimitexperiment}) is equivariant: 
$$
\mbox{for every $u\in\bbr$}\;: \quad \call\left( u^*-u \mid \wt P_u \right) \;\;=\;\;   \call\left( u^* \mid \wt P_0 \right) \;. 
$$
b)~ The sequence of BE  $\vth^*_{nT}$ at stage $n$ defined by (\ref{BE})  is asymptotically as $\nto$ equivariant in the local models $\cale_n^{(\vth)}$ at $\vth$, in the following sense: 
$$
\lim_{\nto}\; \sup_{|u|\le C}\; \left|\;  E_{\vth+\frac{u}{n}} \left( \ell \left( n(\vth^*_{nT}-(\vth{+}\frac{u}{n})) \right)\right) 
\;-\;  E_{\wt P_0} \left( \ell \left( u^* \right)\right)   \;\right| \quad=\quad 0
$$
for every continuous and subconvex loss function $\ell$ which admits a polynomial majorant, and for arbitrary choice of $C<\infty$. \\

An analogous statement also holds for maximum likelhood estimators; we will prove theorem 1.7 in  5.8 below. Now, for quadratic loss, we consider maximal risks over small neighbourhoods of arbitrary points $\vth\in\Theta$, in the sense of contiguous alternatives, and state a local asymptotic minimax theorem (for a particular loss function) which allows to compare arbitrary sequences of $\calg^\xi_{nT}$--measurable estimators $\wt \vth_{nT}$, $n\ge 1$, for the unknown parameter $\vth\in\Theta$.\\

{\bf 1.8 Theorem: } Under Lipschitz and linear growth conditions on $b(\cdot)$ and $\si(\cdot)$, under $(H1)$ and $(H2)$, the following holds for every $\vth\in\Theta$. 

a)~ For squared loss, there is a local asymptotic minimax bound in terms of the BE  $u^*$ of (\ref{BElimitexperiment})  
$$
\lim_{C\uparrow\infty}\; \liminf_{\nto}\; \inf_{\wt \vth_{nT}}\; \sup_{|u|\le C}\; 
E_{\vth+\frac{u}{n}}\left(\; \left[ n\left( \wt \vth_{nT} - (\vth{+}\frac{u}{n}) \right) \right]^2 \;\right) 
\quad\ge\quad E_{\wt P_0}\left( \left[ u^* \right]^2 \right)
$$
where at every stage $n\ge 1$, $\;\inf\limits_{\wt \vth_{nT}}\,$ is with respect to all possible $\calg^\xi_{nT}$--measurable estimators. 

b)~ The BE  sequence $\left(\vth^*_{nT}\right)_n$ of (\ref{BE}) attains the bound given in a), in virtue of 1.7 b).\\

Part a) of this theorem is a result of type [S~85, Cor.\ 62.6], see also [L~72], [M 82] or [V 91]. The proof will be given in 5.9 below. \\

%%%%%%%%%%%%%%%%%%%%%%%%%%%%%%%%%%%%%%%%%%%%%
%\newpage
\section*{2~ Ergodic properties for diffusions with $T$--periodic drift} 
%%%%%%%%%%%%%%%%%%%%%%%%%%%%%%%%%%%%%%%%%%%%%

In this section, we discuss ergodic properties for processes $(\xi_t)_{t\ge 0}$ of type (\ref{process})+(\ref{signal}) under fixed value of the parameter $\vth$. The special features  (\ref{process})+(\ref{signal}) of section 1 are of no importance here. All we assume in the present section is that $(\xi_t)_{t\ge 0}$ is a Markov process with continuous paths, inhomogeneous in time, with  semigroup $(P_{s,t})_{0\le s<t<\infty}$ having measurable densities  
\beqq\label{semigroupdens}
P_{s,t}(x,dy) \;=\; p_{s,t}(x,y)\, dy \quad,\quad 0\le s<t<\infty \;,\; x,y\in\bbr  
\eeqq 
with respect to Lebesgue measure, such that the semigroups is $T$-periodic:  
\beqq\label{semigroupTper}
p_{s,t}(x,y) \;=\; p_{kT+s,kT+t}(x,y)  \quad\mbox{for all $k\in\bbn_0$ and all $0\le s<t<\infty$} \;. 
\eeqq
We use the same notation $ i_T(t) = t \;\mbox{modulo}\; T $ as in (\ref{signal}). Obviously, this setting  covers (\ref{process})+(\ref{signal}) of section 1. The results of section 2 will remain unchanged for Polish state space $(E,\calb(E))$ --~with some reference measure replacing Lebesgue measure in (\ref{semigroupdens})~-- instead of $(\bbr,\calb(\bbr))$. 

We need strong laws of large numbers as $t\to\infty$ for some class of functionals of the process $\xi=(\xi_t)_{t\ge 0}$. This class will be defined in (\ref{property}) below; relevant examples are of the following type.  Let $\Lambda_T(ds)$ denote some $\si$--finite measure on $(-\infty,\infty)$ which is $T$--periodic in the sense that 
\beqq\label{Lambdaper}
\Lambda_T(B)\;=\; \Lambda_T(B+kT)  \quad\mbox{for all $B\in\calb((-\infty,\infty))$ and all $k\in\bbn_0$} \;. 
\eeqq
For $T$--periodic measures $\Lambda_T(ds)$ and for suitable functions $f:\bbr\to\bbr$, define 
\beqq\label{functional-1}
A = (A_t)_{t\ge 0} \quad,\quad    A_t \;=\; \int_0^t f(\xi_s)\, \Lambda_T(ds)  \;,\; t\ge 0  \;. 
\eeqq
%Functionals of this type allow to deal with the periodicity in the process (\ref{process})+(\ref{signal}). 
With fixed $0<r<r'<T$, in view of the log-likelihoods in theorem 1.1,  we may thus consider  
\beqq\label{example-1}
A_t \;=\;  \sum_{k\in\bbn_0 \,,\, kT+r\le t} f( \xi_{kT+r})  
\quad,\quad    \Lambda_T(ds) \;=\; \sum_{k\in\bbn_0} \epsilon_{(kT+r)}(ds)
\eeqq
(with $\epsilon_x$ Dirac measure at the point $x$) or
\beqq\label{example-2}
 A_t \;=\; \int_0^t f(\xi_s)\, 1_{(r,r')}(i_T(s))\, ds 
 \quad,\quad    \Lambda_T(ds) \;=\; \sum_{k\in\bbn_0} 1_{(kT+r,kT+r')}(i_T(s))\; ds  
\eeqq
which both obviously are not additive functionals of the continuous-time process $(\xi_t)_{t\ge 0}$. \\

We fix some notations. We take $\xi$ as  defined on some canonical path space $( \Om, \cala, \bbf, (P_x)_{x\in\bbr})$. 'Almost surely' means almost surely with respect to every $P_x$, $x\in\bbr$. An $\bbf^\xi$--increasing process is an $\bbf^\xi$--adapted c\`adl\`ag process $A = (A_t)_{t\ge 0}$ with nondecreasing paths and $A_0=0$, almost surely. 
We write $(C_T,\calc_T)$ for the space of all continuous functions $\al:[0,T]\to\bbr$ equipped with the metric of uniform convergence and its Borel $\si$-field $\calc_T$. Then $(C_T,\calc_T)$ is a Polish space, and $\calc_T = \si(\pi_t:0\le t\le T)$, the $\si$--field  generated by the coordinate projections $\pi_t$, $0\le t\le T$. \\

The continuous-time Markov process $\xi = (\xi_t)_{t\ge 0}$ induces a $(C_T,\calc_T)$-valued Markov chain  
$$
X = (X_k)_{k\in\bbn_0} \quad,\quad X_k := (\xi_{(k-1)T+s})_{0\le s\le T} \;,\; k\ge 1 \;,\; X_0=\al_0\in C_T   
$$
which we call the chain of  $T$-segments in the path of $\xi$. In virtue of  (\ref{semigroupTper}), the chain $X = (X_k)_{k\in\bbn_0} $ is time homogeneous with one-step-transition kernel  $Q(\cdot,\cdot)$ given  by 
$$
Q(\al,F) \;:=\; P\left(\, (\xi_s)_{0\le s\le T} \in F \mid \xi_0=\al(T)\, \right) \;,\quad \al\in C_T \;,\; F\in\calc_T \;. 
$$
For Harris processes in discrete time, we refer to Revuz [R 75] or Meyn and Tweedie [MT 93]. For Harris processes in continuous time see Az\'ema, Duflo and Revuz [ADR 69] or R\'evuz and Yor [RY~91, Ch.\ 10.3]. Now 'ergodicity' of the process $(\xi_t)_{t\ge 0}$ with $T$-periodic semigroup will be understood as ergodicity of the segment chain $X = (X_k)_{k\in\bbn_0}$: \\

{\bf 2.1 Theorem: } Assume that the embedded chain $(\xi_{kT})_{k\in\bbn_0}$ is positive recurrent in the sense of Harris, and write $\mu$ for its invariant probability on $(\bbr,\calb(\bbr))$.

a) Then the chain $X = (X_k)_{k\in\bbn_0}$  of $T$-segments in the path of $\xi$ is positive recurrent in the sense of Harris. Its  invariant probability  is the unique law $m$ on $(C_T,\calc_T)$ such that 
\beqq\label{invmeasureTper-1}
\left\{ \begin{array}{l}
\mbox{for arbitrary $0=t_0<t_1<\ldots<t_{\ell}<t_{\ell+1}=T$ and $A_i\in\calb(\bbr)$} \;, \\
m\left( \left\{ \pi_{t_i}\in A_i \,,\,0\le i\le \ell{+}1 \right\}\right) \quad\mbox{is given by}\\
\int\ldots\int \mu(dx_0)\, 1_{A_0}(x_0)\, \prod_{i=0}^\ell P_{t_i,t_{i+1}}(x_i,dx_{i+1})\, 1_{A_{i+1}}(x_{i+1}) \;. 
\end{array} \right.
\eeqq

b) For every $\bbf^\xi$--increasing process $A = (A_t)_{t\ge 0}$ with the property
\beqq\label{property}
\left\{ \begin{array}{l}
\mbox{there is some function $F:C_T\to\bbr$, nonnegative, $\calc_T$-measurable, }\\
\mbox{satisfying}\quad   m(F) \;:=\;  \int_{C_T} F\, dm \;<\;  \infty \;, \;\;\mbox{such that} \\
A_{kT} \;=\; \sum\limits_{j=1}^k F(X_k)  
\;=\; \sum\limits_{j=1}^k  F\left(\, (\xi_{(k-1)T+s})_{0\le s\le T} \,\right) \;,\; k\ge 1  
\end{array} \right.
\eeqq
we have the strong law of large numbers 
$$
\lim_{t\to\infty}\, \frac{1}{t}\, A_t \;\;=\;\; \frac{1}{T}\, m(F) \quad\mbox{almost surely} \;. 
$$

\vskip0.8cm
{\bf Proof: } 1) Harris recurrence of the process $(\xi_{kT})_{k\in\bbn_0}$ with invariant probability $\mu$ yields 
$$
A\in\calb(\bbr)\;,\; \mu(A)>0 \quad\mbox{implies}\quad \sum_k 1_A(\xi_{kT}) \;=\;\infty \quad\mbox{almost surely} \;,  
$$
and thus implies that the bivariate chain 
$$
\left(\, \xi_{kT} \,,\, \xi_{(k+1)T} \,\right)_{k\in\bbn_0}
$$
is positive recurrent in the sense of Harris with invariant probability  
$$ 
\mu^{(2)}(dx,dy) \;:=\; \mu(dx)P_{0,T}(x,dy)  \quad\mbox{on $\bbr{\times}\bbr$} \;. 
$$
 
2) Write $m$ for the unique law on $(C_T,\calc_T)$ whose finite dimensional distributions are given by (\ref{invmeasureTper-1}). Since $(C_T,\calc_T)$ is Polish, conditioning with respect to the pair of coordinate projections $(\pi_0,\pi_T)$, the probability $m$ allows for a decomposition  
\beqq\label{invmeasureTper-2}
m(F) \;=\; \int_{\bbr\times\bbr} \mu^{(2)}(dx,dy)\, K( (x,y) , F )
\eeqq
with $K(\cdot,\cdot)$ some transition probability from $\bbr{\times}\bbr$ to $\calc_T$. Comparing with (\ref{invmeasureTper-1}), 
$K((x,y),\cdot)$ is the law of the $\xi$--bridge from $x$ at time $0$ to $y$ at time $T$ 
\beqq\label{invmeasureTper-3}
K( (x,y) , F ) \;:=\; P\left(\, (\xi_s)_{0\le s\le T} \in F \mid \xi_0=x \,,\, \xi_T=y \,\right) \;;  
\eeqq 
for $\mu^{(2)}$--almost all $(x,y)$,  with notations of (\ref{invmeasureTper-1}),  
$\,K\left( (x_0,x_{\ell+1}) , \left\{ \pi_{t_i}\in A_i \,,\,0\le i\le \ell{+}1 \right\}\right)  )$ equals   
$$
\frac{1}{p_{0,T}(x_0,x_{\ell+1})}\, 1_{A_0}(x_0)\, 
\int\ldots\int dx_1\ldots dx_\ell\, \prod_{i=0}^{\ell} p_{t_i,t_{i+1}}(x_i,x_{i+1})\, 1_{A_{i+1}}(x_{i+1})  
$$
whenever $(x_0,x_{\ell+1})$ is in $\{p_{0,T}(\cdot,\cdot)>0\}$, with suitable default definition else. 

3) From (\ref{invmeasureTper-2}) we have for sets $F\in\calc_T$    
\beqq\label{invmeasureTper-4}
m(F)>0 \quad\Lra\quad \left\{ \begin{array}{l} 
\mbox{there is some $\vep>0$ such that}\\
\left\{ (x,y) : K( (x,y) , F )>\vep \right\} \;=:\; B(F,\vep) \\
\mbox{has strictly positive measure under $\mu^{(2)}(dx,dy)$}  
\end{array}\right.
\eeqq
Hence the Harris property of the bivariate chain $\left( \xi_{kT} , \xi_{(k+1)T}\right)_{k\in\bbn_0}$ with invariant measure $\mu^{(2)}(dx,dy)\,$ gives  in combination with  (\ref{invmeasureTper-4}) first 
$$
F\in\calc_T \;,\; m(F)>0 \quad\Lra\quad 
\sum_k 1_{B(F,\vep)}(\xi_{kT},\xi_{(k+1)T}) \;=\;\infty \quad\mbox{almost surely}
$$
and then thanks to  (\ref{invmeasureTper-2})
$$
F\in\calc_T \;,\; m(F)>0 \quad\Lra\quad \sum_k 1_F(X_k) \;=\;\infty \quad\mbox{almost surely} \;. 
$$
We have identified some probability measure $m$ on $(C_T,\calc_T)$ such that sets of positive $m$--measure are visited infinitely often by $X=(X_n)_n$: hence $X$ is Harris. Every  Harris chain admits a unique invariant measure. Periodicity (\ref{semigroupTper}) of the semigroup guarantees that $m$ defined in (\ref{invmeasureTper-1}) is invariant for $X$. 
We have proved that $X=(X_n)_n$ is positive recurrent in the sense of Harris with invariant measure $m$ given by (\ref{invmeasureTper-1}): this is part  a) of the theorem.

5) Consider an $\bbf^\xi$--increasing process $A = (A_t)_{t\ge 0}$ related to a function $F:C_T\to [0,\infty)$  as in (\ref{property}). Then $\Psi = (\Psi_k)_{k\in\bbn_0}$ defined by 
$$
\Psi_k  \;:=\;  A_{kT} \;,\;  k\in\bbn_0
$$ 
is an integrable additive functional of the chain  $X=(X_k)_k$ of $T$--segments in $\xi$. Since $X$ is Harris with invariant measure $m$, we have the ratio limit theorem 
\beqq\label{RLT}
\lim_{k\to\infty}\; \frac{1}{k}\, \Psi_k \;\;=\;\; E_m\left( \Psi_1 \right) \quad\mbox{almost surely} \;. 
\eeqq
But $E_m(\Psi_1)=E_m(F(X_1))=m(F)$, and 
$$
\lim\limits_{k\to\infty}\; \frac{1}{k} \Psi_k   \;\;=\;\;   T\; \lim\limits_{t\to\infty}\; \frac{1}{t}A_t
$$
since $A$ is increasing. This proves part b) of the theorem. \halmos \\

{\bf 2.2 Examples: } Under all assumptions of theorem 2.1,  we deduce in particular  the following laws of large numbers for the functionals  (\ref{functional-1})--(\ref{example-2}). For $0<r<r'<T$ fixed, 

a)~ if $f:\bbr\to\bbr$ is measurable and in $L^1(\, \mu P_{0,r} \,)$, then for $A=(A_t)_t$ considered  in (\ref{example-1})
$$
\lim_{t\to\infty}\;  \frac{1}{t}  \sum_{k\in\bbn_0 \,,\, kT+r\le t} f( \xi_{kT+r})  \;\;=\;\;  \frac{1}{T}\,(\mu P_{0,r})(f)   \quad\mbox{almost surely} \;;
$$
 
b)~if $f:\bbr\to\bbr$ is measurable and in $L^1 \left(\, \int_r^{r'} ds\, (\mu P_{0,s}) \,\right)$, then for the functional  in (\ref{example-2})
$$
\lim_{t\to\infty}\;  \frac{1}{t}  \int_0^t f(\xi_s)\, 1_{(r,r')}(i_T(s))\, ds  \;\;=\;\; \frac{1}{T} \int_r^{r'}  (\mu P_{0,s})(f)\, ds  
\quad\mbox{almost surely}  
$$

c)~ if $f:\bbr\to\bbr$ is measurable and in $L^1 \left(\, \int_0^T \Lambda_T(ds)\, (\mu P_{0,s}) \,\right)$ for some $T$--periodic measure $\Lambda_T$ as defined in (\ref{Lambdaper}), then we have  for the functional $A=(A_t)_t$   in (\ref{functional-1})
$$
\lim_{t\to\infty}\;  \frac{1}{t}  \int_0^t f(\xi_s)\, \Lambda_T(ds)   \;\;=\;\; \frac{1}{T} \int_0^T  \Lambda_T(ds)\, (\mu P_{0,s})(f)\quad\mbox{almost surely}  \;. 
$$
This follows from theorem 2.1: put $F(\al)=f(\al(r))$ for part a),  $F(\al)=\int_r^{r'}f(\al(s))ds$ for part~b), and $F(\al)=\int_0^T f(\al(s)) \Lambda(ds)$ for part c), with $\al\in C_T$. \halmos\\

In order to calculate explicitely the measures occurring in theorem 2.1 or example 2.2, we have to know the semigroup $(P_{s,t})_{0\le s<t<\infty}$ of the time inhomogeneous diffusion $(\xi_t)_{t\ge 0}$. Only in very few cases this can be done explicitely; we mention an Ornstein-Uhlenbeck type example. \\

{\bf 2.3 Example: } With $\si>0$, $\gamma>0$, and some function $S(\cdot)$ which is $T$-periodic and piecewise continuous, consider an Ornstein-Uhlenbeck type diffusion with $T$--periodic drift
$$
d\xi_t  \;=\;  ( S(t)- \gamma\, \xi_t )\,dt  \;+\; \si\, dW_t  \;,\quad t\ge 0 \;. 
$$
The solution with initial value $\xi_0$ is 
$$
\xi_t \;=\; \xi_0\,e^{-\gamma\, t}  \;+\;  \int_0^t e^{-\gamma\, (t-s)} \left(\, S(s)\, ds + \si\, dW_s \right)
$$
and the transition semigroup $( P_{s,t} )_{0\le s<t<\infty}$ of  $\xi$ is 
$$
P_{s,t}(x,\cdot)  \;=\; \caln\left(\;  x\, e^{-\gamma (t-s)} + \int_0^{t-s} e^{-\gamma v} S(t-v)\, dv \;,\; 
e^{-2\gamma (t-s)}\,\frac{e^{2\gamma (t-s)}-1}{2\gamma}\, \si^2  \;\right) \;. 
$$
Extending $S(\cdot)$ to a $T$--periodic function defined on the whole real axis, we define  
$$
M(r) \;=\; \int_0^\infty e^{-\gamma v} S(r-v)\, dv \;=\; \frac{1}{1 - e^{\gamma T}} \int_0^T e^{-\gamma v} S(r-v)\, dv \quad,\quad r\ge 0
$$
which is $T$-periodic. We have weak convergence 
\beqq\label{limit-1}
P_{0,kT}(x,\cdot)  \;\; \lra\;\;  \caln\left(\;  M(0) \;,\; \frac {\si^2}{2\gamma}\;  \right) \;=:\; \mu  \;, \quad k\to\infty 
\eeqq 
for arbitrary  $x\in\bbr$, and similiarly for fixed $0<s<T$ 
$$
P_{0,kT+s}(x,\cdot)  \;\;\lra\;\; \caln\left(\;  M(s) \;,\; \frac {\si^2}{2\gamma} \;\right)  \;, \quad k\to\infty \;. 
$$
It is easy to see that the measure $\mu$ defined in (\ref{limit-1}) is invariant  for the chain $( \xi_{kT} )_{k\in\bbn_0}$, and that sets of positive Lebesgue measure (hence sets of positive $\mu$-measure) are visited infinitely often by $( \xi_{kT} )_{k\in\bbn_0}$,  for every choice of a starting point. Hence $( \xi_{kT} )_{k\in\bbn_0}$ is positive Harris with invariant probability $\mu$. Hence all conditions of theorem 2.1 are satisfied. The last convergence gives 
$$
\mu\, P_{0,s} \;\;=\;\; \caln\left(\;  M(s) \;,\; \frac {\si^2}{2\gamma}\; \right)
\quad,\quad 0<s<T  \;,
$$
thus all measures occcuring in 2.1 or 2.2 are known explicitely. The mapping $s \to \mu P_{0,s}$ describes an 'oscillating stationary regime'. \halmos\\

%%%%%%%%%%%%%%%%%%%%%%%%%%%%%%%%%%%%%%%%%%%%%%%%%
%\newpage
\section*{3~ An exponential inequality} 
%%%%%%%%%%%%%%%%%%%%%%%%%%%%%%%%%%%%%%%%%%%%%%%%%%

In this section, we give an exponential inequality for processes of type (\ref{process})+(\ref{signal}). 
Parameter $\vth$ and periodicity $T$ will play no role in this section, and we consider any time dependent diffusion  
\beqq\label{xiforsection3}
d\xi_t \;=\; b(t,\xi_t)\, dt \;+\; \si(\xi_t)\, dW_t \;,\quad t\ge 0
\eeqq
where $(t,x)\to b(t,x)$ is continuous in restriction to segments $[d_n,d_{n+1}[{\times}\bbr$, for some given deterministic sequence  $(d_n)_n$ with $d_0=0$ and $d_n\uparrow\infty$, under the assumptions 
 \beqq\label{lipschitz}
| b(t,x)-b(t,x') |    +   | \si (x)-\si (x') |    \;\le\; L\, | x-x' | \quad\mbox{for all $x,x'\in\bbr$ and all $t\ge 0$} \;, 
\eeqq
\beqq\label{lingrowth}
| b(t,x) |  \;\le\; L\, ( 1+|x| )  \quad\mbox{for all $x\in\bbr$ and all $t\ge 0$} \;,   
\eeqq
\beqq\label{sigmabounded}
|\si(x)|\;\le\;  M  \quad\mbox{for all $x\in\bbr$} \;.  
\eeqq
Then  a strong solution exists for  (\ref{xiforsection3}) --~we construct it successively on  $[d_n,d_{n+1}]$ as in Karatzas and Shreve [KS 91, theorems 5.2.9+5.2.13], taking for $n\ge 1$ the terminal value of the preceding step as starting value for the following one~-- and is pathwise unique. Clearly this setting includes as a special case equation (\ref{process})+(\ref{signal}) if we put  $b(t,x):=S(\vth,t)+b(x)$, under Lipschitz conditions on $b(\cdot)$ and $\si(\cdot)$, if $\si(\cdot)$ is bounded. 
Up to  a minor modification which allows for time dependence in the drift, the following result is due to Brandt ([B 05], Lemma 2.2.4) who extended  the classical Bernstein inequality for local martingales (see formula (1.5) in Dzhaparidze and van Zanten  [DvZ~01]) to solutions of SDE.\\

{\bf 3.1 Lemma: } (Brandt [B 05]) Fix $0<\la<\frac12$ and $\frac12<\eta<1{-}\la$. 
For $(\xi_t)_{t\ge 0}$ of (\ref{xiforsection3})--(\ref{sigmabounded}) 
%and for $L$ of (\ref{lipschitz})+(\ref{lingrowth}), 
there is some $\Delta_0>0$ (depending only on $L$  of (\ref{lipschitz})+(\ref{lingrowth}), on $\la$ and on $\eta$) such that     
\beqq\label{brandt}\quad
P\left( \sup_{t_1\le t\le {t_1+\Delta}}\, |\xi_t-\xi_{t_1}| > \Delta^\la  \;,\; |\xi_{t_1}|\le   \left(\frac{1}{\Delta}\right)^\eta \right) 
\quad\le\quad
c_1\cdot \exp\left\{-\, c_2\, \left(\frac{1}{\Delta}\right)^{1-2\la} \right\}
\eeqq
holds for  all $0\le t_1<\infty$ and all $0<\Delta<\Delta_0$ , 
with positive constants $c_1$ and $c_2$ which do not depend on $t_1\ge 0$ or on $\Delta\in(0,\Delta_0)$. \\

{\bf Proof: } The proof is from  [B 05].  We  start for $0\le t_1<t<t_1+\Delta$ from   
$$
| \xi_t-\xi_{t_1} |  \;\;\le\;\;  | \int_{t_1}^t \si(\xi_s)\, dW_s | \;+\; \int_{t_1}^t |b(s,\xi_s)|\, ds \;. 
$$
Applying to  $|b(s,\xi_s)| \le |b(s,\xi_{t_1})| + |b(s,\xi_s)-b(s,\xi_{t_1})|$ the  conditions (\ref{lipschitz})+(\ref{lingrowth}), this gives 
$$
| \xi_t-\xi_{t_1} |  \;\;\le\;\;   \left\{ \sup_{t_1\le t\le {t_1+\Delta}}  | \int_{t_1}^t \si(\xi_s)\, dW_s |   \;+\;  L(1+|\xi_{t_1}|)\Delta \right\} \;+\;  L \int_{t_1}^t | \xi_s - \xi_{t_1} |\, ds \;. 
$$
With Gronwall pathwise in $\om$  we obtain (e.g.\ Bass [B 98, lemma I.3.3])
$$
| \xi_t-\xi_{t_1} |  \;\;\le\;\;  \left\{ \sup_{t_1\le t\le {t_1+\Delta}}  | \int_{t_1}^t \si(\xi_s)\, dW_s |   \;+\;  L(1+|\xi_{t_1}|)\Delta \right\}  \cdot e^{Lt} 
$$
 for $0\le t_1<t<t_1+\Delta$. Hence 
\beao
&&P\left(  \sup_{t_1\le t\le {t_1+\Delta}}\, |\xi_t-\xi_{t_1}| > \Delta^\la  \;,\; |\xi_{t_1}|\le   \left(\frac{1}{\Delta}\right)^\eta \right) \\
&&\le\; P\left(  \sup_{t_1\le t\le {t_1+\Delta}}\,   | \int_{t_1}^t \si(\xi_s)\, dW_s | \;>\; \left[ \Delta^\la e^{-L\Delta} - L(1+|\xi_{t_1}|)\Delta  \right]^+ \;,\; |\xi_{t_1}|\le \Delta^{-\eta}   \right) \\
&&\le\; P\left(  \sup_{t_1\le t\le {t_1+\Delta}}\,   | \int_{t_1}^t \si(\xi_s)\, dW_s | \;>\; \left[ \Delta^\la e^{-L\Delta} - L\Delta -  L\Delta^{1-\eta}  \right]^+  \right) 
\eeao
Exploiting assumption (\ref{sigmabounded}), the classical Bernstein inequality for continuous local martingales (formula (1.5) in  [DvZ 01)]) gives
$$
P\left(  \sup_{t_1\le t\le {t_1+\Delta}}\,   | \int_{t_1}^t \si(\xi_s)\, dW_s | \;>\; z  \right)  \;\;\le\;\;  c_1\, e^{-\tilde c_2\, z^2/\Delta } \;,\quad z>0
$$
with positive constants $c_1$ and $\tilde c_2$ which do not depend on $t_1$ or on $\Delta$. We have also 
$$
z \;:=\; \left[ \Delta^\la e^{-L\Delta} - L\Delta -  L\Delta^{1-\eta}  \right]^+  \;=\;  \left[ \Delta^\la ( e^{-L\Delta}-L\Delta^{1-\eta-\la} - L\Delta^{1-\la} )  \right]^+  
\;\ge\quad   \frac12\, \Delta^\la  
$$
provided $\Delta$ is sufficiently small, since  $1-\eta-\la>0$ by assumption. The assertion follows.\halmos\\

%%%%%%%%%%%%%%%%%%%%%%%%%%%%%%%%%%%%%%%%%%%
%
\section*{4~ A result on finite-dimensional convergence}  
%
%      new title for the revised version of this chapter (revision 05.09.09, here 20.10.09)  
%
%%%%%%%%%%%%%%%%%%%%%%%%%%%%%%%%%%%%%%%%%%%%

In this section, we prove two limit theorems (4.1, 4.3, with remark 4.4 below) which will allow to work with log-likelihoods in local models at $\vth$. For fixed value of the parameter $\vth$ which is suppressed from notation, we consider the diffusion $\xi=(\xi_t)_{t\ge 0}$ of (\ref{process})
$$
d\xi_t  \;=\;  \left[ S(t) + b(\xi_t) \right]dt   \;+\;  \si(\xi_t)\,dW_t  
$$
where %--~sligthly more general than (\ref{signal})~-- 
$S(\cdot)$ is some deterministic function with the property  
\beqq\label{Vper}
S  :   [0,\infty)\to\bbr  \quad\mbox{is  $T$-periodic and piecewise continuous} \;. 
\eeqq
This covers (\ref{signal}). 
We write $b(t,x)=S(t)+b(x)$ as in (\ref{xiforsection3}),  assume (\ref{lipschitz})+(\ref{lingrowth}), and strengthen (\ref{sigmabounded}) to 
\beqq\label{sigmadoublebounded}
\frac{1}{M}  \;\le\;   |\si(x)|  \;\le\;  M  \quad\mbox{for all $x\in\bbr$} 
\eeqq
for some large $M$. Systematically,  we combine an appropriate control of fluctuations in the process $(\xi_t)_t$ (by lemma 3.1) with strong laws of large numbers for additive functionals of the chain of $T$-segments in $\xi$ (by theorem 2.1), and consider sequences of auxiliary martingales to which we apply the martingale limit theorem.  \\

{\bf 4.1 Theorem: }  For $\xi$  as above, assume that the embedded chain $(\xi_{kT})_k$ is positive Harris recurrent with invariant probability $\mu$. Then we have for  $0<r<T$ fixed and arbitrary  $h>0$ 
 \beam
 \int_0^{nT} \frac{1}{\si^2(\xi_s)}\, 1_{(r, r+\frac{h}{n})}( i_T(s) )\, ds   &=&   h\, \frac{1}{n} \sum_{j=0}^{n-1}  \frac{1}{\si^2(\xi_{jT+r})}   \;+\; o_{P_x}(1)   \label{thm2-1}\\
  \int_0^{nT} \frac{1}{\si^2(\xi_s)}\, 1_{(r-\frac{h}{n},r)}( i_T(s) )\, ds   &=&   h\, \frac{1}{n} \sum_{j=0}^{n-1}  \frac{1}{\si^2(\xi_{jT+r})}   \;+\; o_{P_x}(1) \label{thm2-2}
 \eeam
as $\nto$, for all $x\in\bbr$, with leading term   
$$
h\, \frac{1}{n} \sum_{j=0}^{n-1}  \frac{1}{\si^2(\xi_{jT+r})}   \quad\lra\quad     h\;  (\mu P_{0,r})(\frac{1}{\si^2})     
$$  
almost surely as $\nto$, as in example 2.2 a). \\

 {\bf Proof: }  The proof of  the approximations (\ref{thm2-1})+(\ref{thm2-2}) is in several steps. We take $n$ large enough to have  $r-\frac{h}{n},r+\frac{h}{n}$ in $(0,T)$, and fix  $0<\la<\frac12$ and $\frac12<\eta<1{-}\la$ as in lemma 3.1. 
 
 1)~ With arbitrary constants $K$, the following auxiliary result will be needed frequently: 
 \beqq\label{fact-1}
 \frac{1}{n} \sum_{j=0}^{n-1} 1_{ \{\, \sup\limits_{0\le s\le T}\, |\xi_{jT+s}| \; >\;  K\, n^\eta  \,\}  }   \;\;\lra\;\;  0   \quad\mbox{almost surely as $\nto$} \;. 
 \eeqq
This is easily derived from the strong law of large numbers in the chain $X=(X_k)_{k\in\bbn_0}$ of $T$-segments in the process $\xi$:   for any fixed $c<\infty$, theorem~2.1~b) with $F(\al)=\sup\limits_{0\le s\le T}|\al(s)|$ (which is a continuous function on $C_T$) gives
$$
 \frac{1}{n} \sum_{j=0}^{n-1} 1_{ \{\, \sup\limits_{0\le s\le T} |\xi_{jT+s}| \; >\;  c \,\}  }    \;\;=\;\;   \frac{1}{n} \sum_{k=0}^{n-1} 1_{ \{F(X_k)>c \}  } 
 \;\;\lra\;\;  m(\{ \al\in C_T :  \sup\limits_{0\le s\le T} |\al(s)| > c \})    
 $$
almost surely as $\nto$. Since the limit on the right hand side decreases to $0$ as $c\uparrow\infty$ and since $n^\eta$ exceeds any fixed level as $n$ tends to $\infty$, (\ref{fact-1}) is proved. 

2)~ Next we prove  for arbitrary $h>0$ and arbitrary starting point $x$ the approximation 
\beqq\label{fact-2}
\frac{1}{n} \sum_{j=0}^{n-1}  \frac{1}{\si^2(\xi_{jT+r - \frac{h}{n}})}  \;\;=\;\; \frac{1}{n} \sum_{j=0}^{n-1}  \frac{1}{\si^2(\xi_{jT+r})} \;+\; o_{P_x}(1) \quad,\quad \nto \;. 
\eeqq
Put $\Delta_n=\frac{h}{n}$. Since by assumption $\si(\cdot)$ is Lipschitz and bounded away from $0$ and $\infty$, the $j$-th summand contributing to the difference in (\ref{fact-2})
 $$
 \left| \frac{1}{\si^2(\xi_{jT+r})} -  \frac{1}{\si^2(\xi_{jT+r - \frac{h}{n}})}  \right|  
 $$
 admits for every $j=0,1,\ldots, n{-}1$ fixed a bound  of type 
 \beao
&&d_1\cdot 1_{ \{\, \sup\limits_{0\le s\le T}\, |\xi_{jT+s}| \;\;>\;\;   \left(\frac{1}{\Delta_n}\right)^\eta  \,\}  }      \\
&&+\quad d_1\cdot 1_{ \{\,   \sup\limits_{jT+r-\frac{h}{n}\le t\le jT+r}\, |\xi_t-\xi_{jT+r-\frac{h}{n}}|  \;\;>\;\; \Delta_n^\la  \; \;,\; \; |\xi_{jT+r-\frac{h}{n}}|  \;\;\le\;\;   \left(\frac{1}{\Delta_n}\right)^\eta     \,\}  }  \\
 &&+\quad d_2\, \Delta_n^\la    \cdot    1_{ \{\,   \sup\limits_{jT+r-\frac{h}{n}\le t\le jT+r}\, |\xi_t-\xi_{jT+r-\frac{h}{n}}| \;\; \le\;\; \Delta_n^\la  \,\}  } 
 \eeao
with suitable constants $d_1$, $d_2$ (where $d_2$ involves the Lipschitz constant $L$).  By the first type of bound combined with step 1, we see that 
$$
 \frac{1}{n} \sum_{j=0}^{n-1}  \left| \frac{1}{\si^2(\xi_{jT+r})} -  \frac{1}{\si^2(\xi_{jT+r - \frac{h}{n}})}  \right| \; 1_{ \{\, \sup\limits_{0\le s\le T}\, |\xi_{jT+s}| \;\;>\;\;   \left(\frac{1}{\Delta_n}\right)^\eta  \,\}  }   
$$
 vanishes almost surely as $\nto$. Next, the exponential inequality in lemma 3.1 (applied to $t_1 = jT+r-\frac{h}{n}$ for $j=0,1,\ldots, n{-}1$) implies that  
 $$
 P\left(\,    \sup\limits_{jT+r-\frac{h}{n}\le t\le jT+r}\, |\xi_t-\xi_{jT+r-\frac{h}{n}}|  \;> \; \Delta_n^\la  \; \;,\; \; |\xi_{jT+r-\frac{h}{n}}|  \;\le \;   \left(\frac{1}{\Delta_n}\right)^\eta  \; ,\;\mbox{some}\; \; j=0,1,\ldots,n{-}1 \,\right)   
 $$
 vanishes as $\nto$. Hence, by the second type of bound,  the probability to find any strictly positive summand in 
 $$
 \frac{1}{n} \sum_{j=0}^{n-1}  \left| \frac{1}{\si^2(\xi_{jT+r})} -  \frac{1}{\si^2(\xi_{jT+r - \frac{h}{n}})}  \right| \; 1_{  \{\,   \sup\limits_{jT+r-\frac{h}{n}\le t\le jT+r}\, |\xi_t-\xi_{jT+r-\frac{h}{n}}|  \;\;>\;\; \Delta_n^\la  \; \;,\; \; |\xi_{jT+r-\frac{h}{n}}|  \;\;\le\;\;   \left(\frac{1}{\Delta_n}\right)^\eta     \,\}   }  
$$
tends to $0$ as $\nto$: hence this sum vanishes in probability as $\nto$.  Finally, by the third type of bounds, we are left to consider averages 
 $$
 \frac{1}{n} \sum_{j=0}^{n-1}  \left| \frac{1}{\si^2(\xi_{jT+r})} -  \frac{1}{\si^2(\xi_{jT+r - \frac{h}{n}})}  \right| \;  1_{ \{\,   \sup\limits_{jT+r-\frac{h}{n}\le t\le jT+r}\, |\xi_t-\xi_{jT+r-\frac{h}{n}}| \;\;\le\;\; \Delta_n^\la  \,\}  } 
 $$
which are bounded by $d_2 \Delta_n^\la$, and thus vanish  as $\nto$. We have proved (\ref{fact-2}). 
 
3)~ Next we show that for arbitrary $h>0$ and arbitrary starting point $x$ 
\beqq\label{fact-3}
\sum_{j=0}^{n-1}  \int_{jT+r-\frac{h}{n}}^{jT+r} \frac{1}{\si^2(\xi_s)}\, ds   
\;\;=\;\; h\, \frac{1}{n} \sum_{j=0}^{n-1}  \frac{1}{\si^2(\xi_{jT+r-\frac{h}{n}})} \;+\; o_{P_x}(1)
\eeqq
as $\nto$. The proof of (\ref{fact-3}) follows the same scheme as the proof of (\ref{fact-2}): we consider for $j=0,1,\ldots, n{-}1$ summands 
\beqq\label{summand-step3}
\int_{jT+r - \frac{h}{n}}^{jT+r} \left| \frac{1}{\si^2(\xi_{s})} -  \frac{1}{\si^2(\xi_{jT+r - \frac{h}{n}})}  \right| ds
\eeqq
and have for these --~since $\si(\cdot)$ is Lipschitz and bounded away from $0$ and $\infty$~--  bounds of type  
 \beao
&&\frac1n\, d_1\cdot 1_{ \{\, \sup\limits_{0\le s\le T}\, |\xi_{jT+s}| \;\;>\;\;   \left(\frac{1}{\Delta_n}\right)^\eta  \,\}  }      \\
&&+\quad \frac1n\, d_1\cdot 1_{ \{\,   \sup\limits_{jT+r-\frac{h}{n}\le t\le jT+r}\, |\xi_t-\xi_{jT+r-\frac{h}{n}}|  \;\;>\;\; \Delta_n^\la  \; \;,\; \; |\xi_{jT+r-\frac{h}{n}}|  \;\;\le\;\;   \left(\frac{1}{\Delta_n}\right)^\eta     \,\}  }  \\
 &&+\quad \frac1n\, d_2\, \Delta_n^\la    \cdot    1_{ \{\,   \sup\limits_{jT+r-\frac{h}{n}\le t\le jT+r}\, |\xi_t-\xi_{jT+r-\frac{h}{n}}| \;\; \le\;\; \Delta_n^\la  \,\}  } 
 \eeao
which allow to proceed  in complete analogy to step 2) above to establish (\ref{fact-3}). 
 
4)~ Combining (\ref{fact-2})+(\ref{fact-3}), we have proved (\ref{thm2-2}). The proof of (\ref{thm2-1}) is similiar. 
%, along the lines of steps 2)+3) above with  obvious notational changes. 
\halmos\\

{\bf 4.2 Corollary: }  For $\xi$  as above, assume that the embedded chain $(\xi_{kT})_k$ is positive Harris recurrent with invariant probability $\mu$. Then we have for  $0<r<T$ fixed and arbitrary  $h>0$ 
\beao
\int_0^{nt} \frac{1}{\si^2(\xi_s)}\, 1_{(r, r+\frac{h}{n})}( i_T(s) )\, ds   &\lra&   \frac{t}{T}\; h\; (\mu P_{0,r})(\frac{1}{\si^2})  \\
\int_0^{nt} \frac{1}{\si^2(\xi_s)}\, 1_{(r-\frac{h}{n},r)}( i_T(s) )\, ds   &\lra&   \frac{t}{T}\; h\;  (\mu P_{0,r})(\frac{1}{\si^2})  
\eeao
almost surely as $\nto$, for every $t>0$ fixed. \\

{\bf Proof: } This is a variant of the preceding proof: for $t>0$ fixed and $\nto$, we consider a sequence $(m_n)_n$ such  that  $m_n\le \frac{tn}{T}< m_n+1$, and copy the proof of  (\ref{thm2-2}) to obtain 
$$
\int_0^{m_n T} \frac{1}{\si^2(\xi_s)}\, 1_{(r-\frac{h}{n},r)}( i_T(s) )\, ds   \quad=\quad   
h\, \frac{1}{n} \sum_{j=0}^{m_n-1}  \frac{1}{\si^2(\xi_{jT+r})}   \;+\; o_{P_x}(1)  
$$
almost surely as $\nto$, for arbitrary $x\in\bbr$. Then thanks to theorem 2.1, the assertion follows as in example 2.2 b). \halmos \\

In the following, we view the process $\xi$ of (\ref{process})+(\ref{Vper}) again as canonical process on the canonical path space $(C,\calc)$, and write $(\eta_t)_{t\ge 0}$ as in the beginning of section 1. 
For $0<r<T$ and $n$ large enough, define martingales $(Y_t^{n,r,h})_{t\ge 0}$ with respect to $\,\bbg^n$ 
$$
\bbg^n = (\calg_{nt})_{t\ge 0} \quad,\quad \calg_t := \bigcap\limits_{\vep>0} \si( \eta_s : 0\le s\le t+\vep)
$$
by 
$$
Y_t^{n,r,h}  \;:=\;  \left\{ \begin{array}{ll}
\int_0^{nt} \frac{1}{\si(\eta_s)}\, 1_{\left( r \,,\, r+\frac{h}{n} \right)}(i_T(s))\, dB_s   & \mbox{if $h>0$} \\
\int_0^{nt} \frac{1}{\si(\eta_s)}\, 1_{\left( r-\frac{|h|}{n} \,,\, r \right)}(i_T(s))\, dB_s   & \mbox{if $h<0$} \;. 
\end{array}\right.   
$$ 
Here $n$ large enough means $\frac{|h|}{n}<\min(r,T-r)$. 
Fix a set of points $0=r_0<r_1<\ldots <r_\ell<r_{\ell+1}=T$ and a set of points ${-}\infty<h_1<h_2<\ldots < h_m<\infty$ in $\bbr\setminus\{0\}$. For $n$ large enough (such that $\frac1n \max\limits_{k=1,\ldots,m}|h_k| < \min\limits_{j=0,\ldots,\ell}(r_{j+1}-r_j)\;$)  compose $\,\bbg^n$-martingales 
\beao
\mathbb{Y}^{n,j} \;:=\; \left( \begin{array}{l} Y_t^{n,r_j,h_1} \\ \ldots \\ Y_t^{n,r_j,h_m} \end{array}\right)_{t\ge 0}
\quad\mbox{$\bbr^m$-valued}  \quad,\quad 
\mathbb{Y}^{n} \;:=\; \left( \begin{array}{l} \mathbb{Y}^{n,1}_t \\ \ldots \\ \mathbb{Y}^{n,\ell}_t \end{array}\right)_{t\ge 0}
\quad\mbox{$\bbr^{\ell m}$-valued} \;. 
\eeao

\vskip0.5cm
{\bf 4.3 Theorem: } a)~For $j=1,\ldots, \ell$, we have weak convergence in $\bbr^m$ as $\nto$ 
$$
\mathbb{Y}^{n,j}_T  \quad\lra\quad  \caln\left( \,0\,,\, \mathbb{A}\, \Gamma_j \,\right) 
\quad,\quad \Gamma_j := (\mu P_{0,r_j})(\frac{1}{\si^2})
$$
where $\mathbb{A}=(A_{i,i'})_{1\le i,i'\le m}$ is the matrix with entries 
$$
A_{i,i'} \;=\; \left\{\begin{array}{ll} 
h_i\wedge h_{i'} & \mbox{if $h_i>0$ and $h_{i'}>0$} \\
|h_i|\wedge |h_{i'}| & \mbox{if $h_i<0$ and $h_{i'}<0$} \\
0 & \mbox{else} \;. 
\end{array}\right.
$$
b)~We have weak convergence in $\bbr^{\ell m}$ as $\nto$ 
$$
\mathbb{Y}^n_T \quad\lra\quad  \caln\left( \,0\,,\, \mathbb{D} \,\right) \quad,\quad 
\mathbb{D} = \left( \begin{array}{lll}  
\mathbb{A}\, \Gamma_1  &  &  0  \\  & \ldots  & \\ 0 & & \mathbb{A}\, \Gamma_\ell
\end{array} \right) \;. 
$$

{\bf Proof: }  By definition of the $\,\bbg^n$-martingales $Y^{n,r_j,h_i}$ or $Y^{n,r_j,h_{i'}}$, we have for $t\ge 0$ 
$$
\langle\, Y^{n,r_j,h_i} \,,\, Y^{n,r_j,h_{i'}} \,\rangle _t    \;=\;  
\int_0^{nt} \frac{1}{\si^2(\eta_s)}\, 1_{(r_j, r_j+\frac{h_i\wedge h_{i'}}{n})}( i_T(s) )\, ds
$$
if $h_i>0$ and $h_{i'}>0$, and 
$$
\langle\, Y^{n,r_j,h_i} \,,\, Y^{n,r_j,h_{i'}} \,\rangle _t    \;=\;  
\int_0^{nt} \frac{1}{\si^2(\eta_s)}\, 1_{( r_j-\frac{|h_i|\wedge |h_{i'}|}{n} ,  r_j )}( i_T(s) )\, ds
$$
if $h_i<0$ and $h_{i'}<0$, and in all other cases 
$$
\langle\, Y^{n,r_j,h_i} \,,\, Y^{n,r_j,h_{i'}} \,\rangle _t    \;\equiv\;  0 \;. 
$$
By corollary 4.2, we have almost surely as $\nto$ for fixed $t$ 
$$
\langle\, Y^{n,r_j,h_i} \,,\, Y^{n,r_j,h_{i'}} \,\rangle _t    \;\lra\;  
\left\{\begin{array}{ll} 
\frac{t}{T}\; (h_i\wedge h_{i'})\; \Gamma_j  & \mbox{if $h_i>0$ and $h_{i'}>0$} \\
\frac{t}{T}\; (|h_i|\wedge |h_{i'}|)\; \Gamma_j  & \mbox{if $h_i<0$ and $h_{i'}<0$} \\
0 & \mbox{else} \;. 
\end{array}\right.
$$
with $\Gamma_j$ as above. In virtue of the martingale convergence theorem (Jacod and Shiryaev 1987, VIII.3.22), this implies weak convergence 
$$ 
\left( \mathbb{Y}^{n,j}_t \right)_{t\ge 0}  \quad\lra\quad  \left( \frac{1}{T}\, \mathbb{A}\, \Gamma_j\right)^{1/2}\, \mathbb{B}^j
$$
in the Skorohod path space of c\`adl\`ag functions $[0,\infty)\to\bbr^m$, where $\mathbb{B}^j$ is $m$-dimensional Brownian motion and $\mathbb{A}$ as defined above. This implies a). For $j\neq j'$ chose  $\mathbb{B}^j$, $\mathbb{B}^{j'}$ independent, and compose  
$$
\mathbb{B} \;\;:=\;\; \left( \begin{array}{l} \mathbb{B}^1_t \\ \ldots \\ \mathbb{B}^m_t \end{array}\right) _{t\ge 0}
\quad\mbox{$\bbr^{\ell m}$-valued} \;. 
$$ 
Then the martingale convergence theorem gives weak convergence as $\nto$ 
$$
\left(\, \mathbb{Y}^n_t \,\right)_{t\ge 0}  \quad\lra\quad  \left( \frac{1}{T}\, \mathbb{D} \right)^{1/2}\, \mathbb{B} 
$$
in the Skorohod path space of c\`adl\`ag functions $[0,\infty)\to\bbr^{\ell m}$, with $\mathbb{D}$ as defined above, since 
$$
\langle\, Y^{n,r_j,h_i} \,,\, Y^{n,r_{j'},h_{i'}} \,\rangle _t   \;\equiv\; 0 \quad\mbox{whenever $j\neq j'$}  \;.  
$$
This implies assertion b).\halmos \\

{\bf 4.4 Remark: } The matrix $\mathbb{A}$ defined in theorem 4.3 has the structure of the covariance kernel of two-sided one-dimensional Brownian motion. Hence we shall read theorem 4.3 as a result on processes which for fixed $T$ are now indexed in the parameter $h$  
\beqq\label{interpretation-1}
\left( Y^{n,r_1,h}_T \right)_{h\in[-K,K]} \quad, \quad\ldots\quad ,\quad \left( Y^{n,r_\ell,h}_T \right)_{h\in[-K,K]} 
\eeqq
and which converge by theorem 4.3 as $\nto$ in the sense of finite dimensional distributions to 
\beqq\label{interpretation-2}
\Gamma_1^{1/2}\, (\wt W^1_h)_{h\in[-K,K]}  \quad, \quad\ldots\quad ,\quad \Gamma_\ell^{1/2}\, (\wt W^\ell_h)_{h\in[-K,K]} 
\eeqq
with independent two-sided one-dimensional  Brownian motions  $\wt W^j=(\wt W^j_h)_{h\in\bbr}$, $1\le j\le \ell$. Here the choice of determined  intervals $[-K,K]$ is in fact irrelevant. We shall make use of finite dimensional convergence  (\ref{interpretation-1})+(\ref{interpretation-2}) when considering martingale terms in the log-likelihoods of local models at some fixed point $\vth\in\Theta$. It is through the reinterpretation  (\ref{interpretation-1})+(\ref{interpretation-2}) of theorem 4.3  that the local parameter $h$ in local models at $\vth$ begins to play the role of 'time'. \halmos\\

%%%%%%%%%%%%%%%%%%%%%%%%%%%%%%%%%%%%%%%%%%%
%\newpage
\section*{5~ Proofs for the statistical results stated in section 1} 
% 
%  veraendert ab 26.10.09 
%%%%%%%%%%%%%%%%%%%%%%%%%%%%%%%%%%%%%%%%%%%%

In order to prove the results stated in section 1, we consider the process $\xi$ defined by (\ref{process})+(\ref{signal})  which depends on the parameter $\vth\in\Theta$. We assume the Harris condition $(H1)$ for all values of $\vth\in\Theta$, and Lipschitz and linear growth conditions on $b(\cdot)$ and $\si(\cdot)$ in combination with $(H2)$. This implies that  for all $\vth$, all assumptions made in section 2 hold, as well as those of section 3 and of section 4. 
%(\ref{lipschitz})+(\ref{lingrowth}) of section 3 and (\ref{sigmadoublebounded}) of section 4. 
%Whenever we consider sequences of local models at $\vth$, we will need the function $\la^*(\cdot)$ of (\ref{signal}) only on some fixed open set in $[0,T]$: hence we can take $\la^*(\cdot)$ as strictly positive, continuous and $T$-periodic on $\bbr$. 
As in section 1, $\,Q^\vth$ is the law of $\xi$ under $\vth$ on the canonical path space $(C,\calc)$, $\;\eta\,$ the canonical process on $(C,\calc)$, and $\,\bbg\,$  the canonical filtration on $(C,\calc)$. 
%$\bbg = \left(\calg_t\right)_{t\ge 0}$ with $\calg_t = \bigcap_{t'>t}  \si \left\{ \eta_s : 0 \le \eta_s\le t' \right \}$. 
%Whenever convenient, we use default definitions  $Q^\vth:=Q^0$ for $\vth\le 0$ and $Q^\vth:=Q^{T-a}$ for $\vth\ge T{-}a$, in obvious extension of (\ref{signal}) with $\Theta=(0,T{-}a)$. 
\\ 

We restart from the representation of likelihood ratios  $L^{\zeta ' / \zeta}_t$ of $\,Q^{\zeta '}|\calg_t$ with respect to $\,Q^\zeta|\calg_t$   given before theorem 1.1: for $\zeta\in\Theta$ and  $\zeta '$ close to $\zeta$, $\;L^{\zeta ' / \zeta}_t$ equals 
\beam\label{LRglobalmodel+}
&&  \exp \left\{\;  
\left[  -\;  \int_0^t \frac{\la^*(s)}{\si(\eta_s)} \, 1_{(\zeta,\zeta ')}(i_T(s))\, dB_s   \;-\; \frac12  \int_0^t   \left(\frac{\la^*(s))}{\si(\eta_s)}\right)^2\, 1_{(\zeta,\zeta ')}(i_T(s))\, ds \right] \right. \\
&&+\; \left. \left[ \int_0^t \frac{\la^*(s)}{\si(\eta_s)} \, 1_{(\zeta+a,\zeta '+a)}(i_T(s))\, dB_s   \;-\; \frac12  \int_0^t   \left(\frac{\la^*(s))}{\si(\eta_s)}\right)^2\, 1_{(\zeta+a,\zeta '+a)}(i_T(s))\, ds \right] \nonumber
\; \right\}
\eeam
in case $\zeta<\zeta '<\zeta+a$, with $B$ Brownian motion under the reference point $\zeta$, and 
\beam\label{LRglobalmodel-}
&&\exp \left\{\;  
\left[ \int_0^t \frac{\la^*(s)}{\si(\eta_s)} \, 1_{(\zeta',\zeta )}(i_T(s))\, dB_s   \;-\; \frac12  \int_0^t   \left(\frac{\la^*(s))}{\si(\eta_s)}\right)^2\, 1_{(\zeta',\zeta)}(i_T(s))\, ds \right] \right. \\
&&+\; \left. \left[ -\; \int_0^t \frac{\la^*(s)}{\si(\eta_s)} \, 1_{(\zeta'+a,\zeta +a)}(i_T(s))\, dB_s   \;-\; \frac12  \int_0^t   \left(\frac{\la^*(s))}{\si(\eta_s)}\right)^2\, 1_{(\zeta'+a,\zeta+a)}(i_T(s))\, ds \right] \nonumber
 \; \right\}   
\eeam
in case $\zeta'<\zeta<\zeta'+a$. We state a lemma related to the geometry --~in Hellinger sense~-- of the model determined by (\ref{LRglobalmodel+})+(\ref{LRglobalmodel-}).   \\

{\bf 5.1 Lemma: }  Consider $\zeta , \zeta'\in\Theta$  such that  $\zeta<\zeta '<\zeta+a$ or $\zeta'<\zeta<\zeta'+a$ holds. Then there are  positive constants $c_j$ and $k$ not depending   on $\zeta , \zeta',t$   such that  for all  $t\ge 0$
\beqq\label{maxbound}
E_\zeta \left(\,  \sup_{0\le s\le t}\,  L_t^{\zeta' / \zeta}  \right)   \quad\le\quad   c_0 \;+\;  c_1\left(\left\lfloor \frac{t}{T}\right \rfloor + 1\right)  |\zeta'-\zeta|     \;, 
\eeqq
\beqq\label{hellinger-1} 
E_\zeta \left(\left[\, 1 \,-\,  \left( L_t^{\zeta' / \zeta}  \right)^{1/2} \right]^2 \right)    \quad\le\quad    \sum_{j=1}^3 \;c_j  \left(\left\lfloor \frac{t}{T}\right \rfloor + 1\right)^j  |\zeta'-\zeta|^j      \;, 
\eeqq
\beqq\label{hellinger-2} 
E_\zeta \left(\left[\, 1 \,-\,  \left( L_t^{\zeta' / \zeta}  \right)^{1/4} \right]^4 \right)    \quad\le\quad    \sum_{j=2}^5 \;c_j  \left(\left\lfloor \frac{t}{T}\right \rfloor + 1\right)^j  |\zeta'-\zeta|^j        \;, 
\eeqq
\beqq\label{hellinger-3} 
E_\zeta \left( \left[  L_t^{\zeta' / \zeta}  \right]^{1/2} \right)  \quad\le\quad   \exp\left\{\, -\, k \left\lfloor \frac{t}{T}\right \rfloor   |\zeta'-\zeta|  \right\}  \;. 
\eeqq
\vskip0.8cm

{\bf Proof: } 1)~ For $\zeta$, $\zeta'$ in $\Theta$, the likelihood ratio process $(L_t^{\zeta' / \zeta})_{t\ge 0}$ under $Q^\zeta$ is the exponential 
\beqq\label{fact-11}
L^{\zeta' / \zeta}   \;\;=\;\;  
\cale_\zeta \left( \int _0^\cdot \delta_s\, dB_s \right)   \quad\mbox{with}\quad \delta_s \;=\;  \frac{S(\zeta ',s)-S(\zeta,s)}{\si(\eta_s)} \;. 
\eeqq
Hence $(L_t^{\zeta' / \zeta})_{t\ge 0}$ under $Q^\zeta$ solves the SDE 
$$
L_t^{\zeta' / \zeta} \;=\; 1 \;+\;  \int _0^t L_s^{\zeta' / \zeta}\; \delta_s\, dB_s \;,\quad t\ge 0 
$$
from which one obtains by Ito formula under $Q^\zeta$ 
\beqq\label{fact-12}
V_t := \left(L_t^{\zeta' / \zeta}\right)^{1/4} \quad\mbox{satisfies}\quad  V_t \;=\; 1 - \frac{3}{32}\int_0^t V_s\, \delta^2_s\, ds + \frac{1}{4}\int_0^t V_s\, \delta_s\, dB_s \;, 
\eeqq
 \beqq\label{fact-13}
V_t := \left(L_t^{\zeta' / \zeta}\right)^{1/2} \quad\mbox{satisfies}\quad  V_t \;=\; 1 - \frac{1}{8}\int_0^t V_s\, \delta^2_s\, ds + \frac{1}{2}\int_0^t V_s\, \delta_s\, dB_s \;. 
\eeqq
 
2)~  From now on in this proof, we fix $\zeta\in\Theta$ and consider $\zeta'$  such that either $\zeta<\zeta '<\zeta+a$ or $\zeta'<\zeta<\zeta'+a$ holds. In the first case we have  according to (\ref{LRglobalmodel+})   
$$
\delta_s \;:=\; \frac{\la^*(s)}{\si(\eta_s)}\, \left[ -\, 1_{ ( \zeta ,  \zeta' )} \,+\,  1_{ ( \zeta+a , \zeta'+a )}  \right](i_T(s))  
$$
with non-overlapping intervals, in the second case  by (\ref{LRglobalmodel-}) 
$$
\delta_s \;:=\; \frac{\la^*(s)}{\si(\eta_s)}\, \left[\, 1_{ ( \zeta' ,  \zeta )} \,-\,  1_{ ( \zeta'+a , \zeta+a )}  \right](i_T(s)) \;. 
$$
Let $s\to d(s)$ denote the deterministic $T$--periodic function 
$$
d(s) \;:=\; 1_{\{ \delta \neq 0 \}}(s)  \;=\;  \left[\, 1_{ ( \zeta' ,  \zeta )} +  1_{ ( \zeta'+a , \zeta+a )}  \right](i_T(s))
\;,\; s\ge 0 \;. 
$$
By our assumptions  on $\la^*(\cdot)$  and $\si(\cdot)$, the process $(\frac{\la^*(s)}{\si(\eta_s)})_{s\ge 0}$ is bounded  away from both $0$ and $\infty$. Hence there are some $0<\ul{c}<\ol{c}<\infty$  not depending   on $\zeta , \zeta' , t$ such that 
\beqq\label{integraldeltapower}
 \ul{c}\, d(s) \;\le\; \delta^2_s  \;\le\; \ol{c}\, d(s)      \quad,\quad  
 2  \left\lfloor \frac{t}{T}\right \rfloor  |\zeta'-\zeta|   \;\le\;  
 \int_0^t d(s)\, ds \;\le\;   2  \left(\left\lfloor \frac{t}{T}\right \rfloor + 1\right)  |\zeta'-\zeta|   
\eeqq
for all $t>0$; here $\lfloor x \rfloor$ denotes the biggest integer stricly smaller than~$x$.

3)~ We  prove (\ref{maxbound}). The process $V_t := \left(L_t^{\zeta' / \zeta}\right)^{1/2}$ of (\ref{fact-13}) is nonnegative and satisfies  $E_\zeta(V^2_s)=1$ for all $s$. Thus  (\ref{fact-13}) gives 
$$
0  \;\le\;  V_t   \;\le\; 1 \;+\; \frac12 \int_0^t V_s\, \delta_s\, dB_s
$$
and allows to write (with [IW 89, p.\ 110])
\beao
E_\zeta \left(\,  \sup_{0\le s\le t}\,  V_t^2  \,\right)    
&\le&  c\left(\, 1 \,+\,  E_\zeta \left(\,   \sup_{0\le s\le t}\,  \left(   \int_0^s V_r\, \delta_r\, dB_r  \right)^2 \right) \right) \\
&\le&  c' \left(\, 1 \,+\,  E_\zeta \left(\,   \int_0^t V^2_s\, \delta^2_s\, ds  \right) \right) \\
&\le&  c'' \left(\, 1 \,+\,  E_\zeta \left(\,   \int_0^t V^2_s\, d(s)\, ds  \right) \right)  
\quad\le\quad    \wt c \left(\, 1 \,+\, \int_0^t d(s)\, ds \right)  
\eeao
which combined with (\ref{integraldeltapower}) gives (\ref{maxbound}). 

4)~ We  prove (\ref{hellinger-1}). We start from $V_t := \left(L_t^{\zeta' / \zeta}\right)^{1/2}$ for which  (\ref{fact-13}) yields the bound 
\beqq\label{bound-10}
 -\, \frac12  \int_0^t V_s\, \delta_s\, dB_s    \;\;\le\;\;  1-V_t   \;\le\;  \frac{1}{8}\int_0^t V_s\, \delta^2_s\, ds \;+\; \left(\frac{1}{2}\int_0^t V_s\, \delta_s\, dB_s\right)^-
\eeqq
where $f^-$ denotes the negative part in $f=f^+-f^-$. Combining (\ref{maxbound})+(\ref{integraldeltapower}) we have 
\beao
E_\zeta \left( \left[ \int_0^t V_s\, \delta^2_s\, ds \right]^2 \right)  &\le&  E_\zeta \left( \sup_{0\le s\le t}\, V^2_s \right)\;  \left[  \ol{c}  \int_0^t d(s)\, ds \right]^2  \\
&\le&  c_2 \left(\left\lfloor \frac{t}{T}\right \rfloor + 1\right)^2  |\zeta'-\zeta|^2 \;+\;  c_3 \left(\left\lfloor \frac{t}{T}\right \rfloor + 1\right)^3  |\zeta'-\zeta|^3
\eeao
for suitable $c_2,c_3$. For  the martingale term in (\ref{bound-10}),  using again  (\ref{integraldeltapower}) and $E_\zeta(V^2_s)=1$, we have 
$$
E_\zeta \left( \left[ \int_0^t V_s\, \delta_s\, dB_s \right]^2 \right)   \;=\;  E_\zeta \left( \int_0^t V^2_s\, \delta^2_s\, ds \right)  \;\le\; 
\ol{c}\, \int_0^t d(s)\, ds  
\;\;\le\;\;  c_1  \left(\left\lfloor \frac{t}{T}\right \rfloor + 1\right)  |\zeta'-\zeta|
$$
for all $t\ge 0$.  Squaring  the bound (\ref{bound-10}), both inqualities together give (\ref{hellinger-1}).

5)~ To  prove (\ref{hellinger-2}), we define $V_t := \left(L_t^{\zeta' / \zeta}\right)^{1/4}$ and have from  (\ref{fact-12})  up to different constants again a bound of form  (\ref{bound-10}) . We use 
$$
E_\zeta \left( \left[ \int_0^t V_s\, \delta^2_s\, ds \right]^4 \right)  
\;\;\le\;\;  E_\zeta \left( \sup_{0\le s\le t}\, V^4_s \right)  \left[ \ol{c} \int_0^t d(s)\, ds \right]^4  
$$
and with [IW 89, p.\ 110]
\beao
E_\zeta \left( \left[ \int_0^t V_s\, \delta_s\, dB_s \right]^4 \right)  &\le& c\, E_\zeta \left(  \left[ \int_0^t V^2_s\, \delta^2_s\, ds \right]^2   \right)
\;\;\le\;\;  c\, E_\zeta \left(  \sup_{0\le s\le t}\, V^4_s \right) \left[ \ol{c} \int_0^t d(s)\, ds \right]^2    \;. 
\eeao
Applying (\ref{maxbound}) to the first factor and (\ref{integraldeltapower}) to the second factor appearing on the right hand side,  the sum of the right hand sides of both last inequalities is at most 
$$
\sum_{j=2}^5  \;c_j \left(\left\lfloor \frac{t}{T}\right \rfloor + 1\right)^j  |\zeta'-\zeta|^j \;. 
$$

6)~ To prove (\ref{hellinger-3}), we deduce  from the exponential  representation  (\ref{fact-11})  under $Q^\zeta$ that 
$$
\left(  L_t^{\zeta' / \zeta}  \right)^{1/2}    \;\;=\;\;    
\cale_\zeta \left( \int _0^\cdot (\frac12 \delta_s)\, dB_s \right)_t \; \exp\left\{\, -\frac18 \int_0^t \delta^2_s\, ds\,\right\} \;. 
$$
So  the  lower bound in (\ref{integraldeltapower}) applies and shows 
$$
E_\zeta \left(  \left(  L_t^{\zeta' / \zeta}  \right)^{1/2}  \right)  \;\;\le\;\;  \exp\left\{\, -\frac18\; \ul{c} \int_0^t d(s)\, ds\,\right\} 
 \;\;\le\;\;  \exp\left\{\, -\, k \left\lfloor \frac{t}{T}\right \rfloor   |\zeta'-\zeta|  \right\} 
$$
for suitable $k>0$. We have proved all assertions of lemma 5.1.\halmos\\

From (\ref{hellinger-1}) in lemma 5.1, we have bounds for the  squared  Hellinger distance  when we observe the  trajectory of  $(\xi)_{t\ge 0}$  up to time ${nT}$: for $\zeta'$ sufficiently close to $\zeta$, 
\beqq\label{hellinger-4}
H^2\left(\, Q^\zeta | \calg_{nT} \,,\, Q^{\zeta'} | \calg_{nT} \,\right)  \;=\;   \frac12\, E_\zeta \left( \left[ 1 - \left(L_{nT}^{\zeta' / \zeta}\right)^{1/2} \right]^2 \right)    \;\;\le\;\;   C\, \sum_{j=1}^3  n^j\,  |\zeta'-\zeta|^j 
\eeqq
for all $n\ge 1$, with some constant $C$ which does not depend on  $\zeta,\zeta',n$. In particular, the parametrization in $\left\{ Q^\zeta | \calg_{nT} : \zeta\in\Theta \right\}$ is --~in Hellinger distance~-- H\"older continuous with index $\frac12$.  In local models at some fixed reference point in the parameter space, the right hand side of (\ref{hellinger-4}) forces local scale to be  proportional to~$\frac1n$. This leads to local models  $\cale_n^{(\vth)}$ at $\vth$ as considered in 1.1.  \\

{\bf 5.2 Proof of theorem 1.1: } a)~ From  (\ref{LRglobalmodel+})+(\ref{LRglobalmodel-})   we have for  $n$ sufficiently large  
\beam\label{LRlocalmodel+}
L_{nT}^{ (\vth+\frac{h}{n}) / \vth } &=& \exp \left\{\;   
\int_0^{nT} \frac{\la^*(s)}{\si(\eta_s)} \, \left[ 1_{(\vth+a,\vth+a+\frac{h}{n})} - 1_{(\vth,\vth+\frac{h}{n})} \right] (i_T(s))\,  dB_s \right. \\
&&-\; \left. \frac12\, \int_0^{nT}  \left(\frac{\la^*(s))}{\si(\eta_s)}\right)^2\, \left[ 1_{(\vth+a,\vth+a+\frac{h}{n})} + 1_{(\vth,\vth+\frac{h}{n})} \right] (i_T(s))\, ds \nonumber
 \; \right\}   
\eeam
 in case $u=h>0$, where all intervals are disjoint, and in case $u=-h<0$  
\beam\label{LRlocalmodel-}
L_{nT}^{ (\vth-\frac{h}{n}) / \vth } &=& \exp \left\{\;   
\int_0^{nT} \frac{\la^*(s)}{\si(\eta_s)} \, \left[ 1_{(\vth-\frac{h}{n},\vth)} - 1_{(\vth+a-\frac{h}{n},\vth+a)} \right] (i_T(s))\,  dB_s \right. \\
&&-\; \left. \frac12\, \int_0^{nT}  \left(\frac{\la^*(s))}{\si(\eta_s)}\right)^2\, \left[ 1_{(\vth+a-\frac{h}{n},\vth+a)} + 1_{(\vth-\frac{h}{n},\vth)} \right] (i_T(s))\, ds \nonumber
 \; \right\}   \;. 
\eeam
Since the function $\la^*(\cdot)$ is deterministic, $T$-periodic and continuous, since the intervals occurring in the integrands are of length $O(\frac{1}{n})$, the logarithm of both expressions (\ref{LRlocalmodel+})+(\ref{LRlocalmodel-}) behaves as $\nto$ as  
\beqq\label{neuedarstellung}
M^{n,\vth,h}_T \;-\; \frac12 \langle M^{n,\vth,h} \rangle_T \quad,\quad 
M^{n,\vth,h}_t \;:=\; {\rm sgn} (h) \left( \la^*(\vth+a) Y_t^{n,\vth+a,h}  - \la^*(\vth) Y_t^{n,\vth,h} \right) 
\eeqq 
with notations of theorem 4.3, for every fixed value of $h$. Finite-dimensional convergence as $\nto$ 
$$
\left( M^{n,\vth,h}_T \right)_{h\in [-K,K]} \quad\lra\quad 
J_\vth^{1/2}\cdot \left( \wt W_h \right)_{h\in [-K,K]} \;=\; \left(\, \wt W( u J_\vth ) \right)_{u\in [-K,K]}
$$
follows from theorem 4.3 and remark 4.4, with $(\wt W_u)_{u\in\bbr}$ two-sided standard Brownian motion and   
$$
J_\vth  \;=\;  \left\{  (\la^*(\vth))^2\, (\mu^{(\vth)} P^{(\vth)}_{0,\vth})  + ( \la^*(\vth{+}a))^2\,  (\mu^{(\vth)} P^{(\vth)}_{0,\vth+a}) \right\} (\frac{1}{\si^2})   
$$
as in theorem 1.1. With similiar arguments, convergence of angle bracketts in (\ref{neuedarstellung}) is immediate  from theorem 4.1. The actual choice of $[-K,K]$ above being irrelevant for finite-dimensional convergence, we have proved convergence as $\nto$  of the finite dimensional distributions for likelihood ratios in local models at $\vth$ 
$$
\left(\; L_{nT}^{(\vth + \frac{u}{n}) / \vth } \;\right) _{u \in \Theta_{\vth,n}}  \;\mbox{under $Q^\vth$}
\quad,\quad  \Theta_{\vth,n} := \{ u\in\bbr : \vth + \frac{u}{n} \in\Theta\}  
$$
to the finite dimensional distributions of likelihoodratios (\ref{limitlikelihoodratio})+(\ref{limitexperiment}) in the limit model $\wt\cale$
$$
\left(\wt L^{ u / 0 }\right)_ {u\in\bbr}  \;\mbox{under $\wt P_0$}
\quad,\quad   \wt L^{ u / 0 } \;:=\;  \exp \left\{\;   \wt W(u J_\vth)  -  \frac{1}{2} | u J_\vth |  \;\right\} 
$$
as asserted in theorem 1.1 a). 

b)~ Convergence of local models $\cale_n^{(\vth)}$ in the sense of [S 85, p.\ 302] as $\nto$ to the limit model~$\wt\cale$  reduces to finite dimensional convergence of likelihoods  as proved in a) since all probability measures in $\left\{ Q^{\vth+\frac{h}{n}}{\mid}\calf_{nT} : h\in\Theta_{\vth,n} \right\}$ and in $\left\{ \wt P_h : h\in\bbr \right\}$ are equivalent. \halmos\\

We remark that the same scheme of proof based on theorem 4.3  allows to consider signals of type (\ref{signal}) having finitely many discontinuities which arise at epochs $0<\vth_1<\vth_2<\ldots<\vth_d<T$, cf.\ remark 1.3 for $d=2$, and leads to $d$-dimensional parameter $\vth=(\vth_1,\ldots,\vth_d)$ and limits of local models at $\vth$ which are $d$-fold products --~with suitable scaling in the factor models~-- of independent experiments $\wt \cale$. \\

%%%%%%%%
% diese zwei lemmata veraendert 05.11.09
%%%%%%%%%
The next two lemmata deal with the limit experiment $\wt \cale$ in (\ref{limitexperiment})+(\ref{limitlikelihoodratio}).\\

{\bf 5.3 Lemma: } In $\wt \cale$, for every $u\in\bbr$ fixed, there is a two-sided  $\wt P_u$--Brownian motion $\wt{\wt W}$ such that 
\beqq\label{equivariance-1}
\wt L^{(u+h) / u } \;=\; \exp\left\{\; \wt{\wt W}(h J_\vth)  - \frac12\, |hJ_\vth| \right\} \quad\mbox{for all $h\in\bbr$}   \;. 
\eeqq
Hence for every  $u\in\bbr$, $\left\{ \wt P_{u+h} : h\in\bbr \right\}$ is statistically the same experiment as $\wt\cale$. \\

{\bf Proof: } We will consider the likelihood ratio $\wt L^{u / u_0}$ separately in all possible cases $0<u_0<u$, $\;0<u<u_0$, $\;u_0<0<u$, $\;u<0<u_0$, $\ldots$. We give the detailed proof in case $u<0<u_0$, using notation (\ref{interpretationPu}) of remark 1.2, and suppressing the scaling factor as there.  Write $(\eta^{(1)},\eta^{(2)})$ for the canonical process on $C([0,\infty),\bbr^2)$. By (\ref{interpretationPu}), the canonical process has semimartingale characteristics 
\beao
&&\left( 
\left(   \begin{array}{l} 0 \\  t\wedge|u| \end{array}   \right) , \left( \begin{array}{ll} t & 0 \\ 0 & t \end{array} \right) , \left(  \begin{array}{l} 0 \\  0 \end{array}  \right)
\right) \quad\mbox{under  $\wt P_u$ with $u<0$} \\
&&\left( 
\left(   \begin{array}{l}  t\wedge u_0 \\ 0 \end{array}   \right) , \left( \begin{array}{ll} t & 0 \\ 0 & t \end{array} \right) , \left(  \begin{array}{l} 0 \\  0 \end{array}  \right)
\right) \quad\mbox{under  $\wt P_{u_0}$ with $u_0>0$}
\eeao
hence the likelihood ratio process of $\wt P_u$ to $\wt P_{u_0}$ relative to the canonical filtration on the path space $C([0,\infty),\bbr^2)$ is given by 
$$
\exp\left\{  \int_0^t ({-}1_{(0,u_0)}(t))\,d\wt m^{(1)}_t  \;+\; \int_0^t 1_{(0,|u|)}(t)\,d\wt m^{(2)}_t 
\;-\; \frac12  (u_0\wedge t) \;-\; \frac12  (|u|\wedge t)   \right\}  \;\;,\;\;  t\ge 0 
$$
where $\wt m^{(i)}$ denotes the martingale part of $\eta^{(i)}$ under $\wt P_{u_0}$.  $\;\wt m^{(1)}$ and $\wt m^{(2)}$  being independent $\wt P_{u_0}$--Brownian motions by (\ref{interpretationPu}),  we form a two-sided $\wt P_{u_0}$--Brownian motion $\wt M$ using $\wt m^{(1)}$ on the right and $\wt m^{(2)}$ on the left branch.  
Letting $t\to\infty$ in the above likelihood, we get 
$$
\wt L^{u / u_0}     \;=\;   \exp\left\{  \left( \wt m^{(2)}_{|u|} - \wt m^{(1)}_{u_0} \right)  \;-\; \frac12 ( |u|+u_0 ) \right\}   \;=\;   
\exp\left\{  \left( \wt M_u - \wt M_{u_0} \right)  \;-\; \frac12 \left| u-u_0 \right| \right\}     
$$
since $u<0<u_0$. Also $\wt{\wt W}$ defined by $\wt{\wt W}(h) :=  \wt M_{u_0+h} - \wt M_{u_0}$, $h\in\bbr$,  is a two-sided $\wt P_{u_0}$--Brownian motion, and the last equation takes the form 
$$
\wt L^{u / u_0}     \;=\;  \exp\left\{   \wt{\wt W}(u-u_0)  \;-\; \frac12 \left| u-u_0 \right| \right\}   \;,\quad u<0<u_0 \;. 
$$
We have proved (\ref{equivariance-1})  in case $u<0<u_0$. The remaining cases are proved similiarly.  
Reintroducing the scaling factor $J_\vth$ of (\ref{limitexperiment})+(\ref{limitlikelihoodratio}) amounts  to replace $u$ and $u_0$ by $u\,J_\vth$ and $u_0\,J_\vth$. \halmos\\

{\bf 5.4 Lemma: } In $\wt \cale$,  
%of (\ref{limitexperiment})+(\ref{limitlikelihoodratio}) for local models at $\vth$, 
we have  for all $u,u'$ in $\bbr$  
\beqq\label{hellinger-11}
E_u\left( \left[ 1 - \left( \wt L^{u' / u} \right)^{1/2} \right]^2 \right) \;\;\le\;\;  \frac14\, J_\vth\, |u'-u|   
\eeqq
\beqq\label{hellinger-12}
E_u\left( \left[ 1 - \left( \wt L^{u' / u} \right)^{1/4} \right]^4 \right) \;\;\le\;\;  c_\vth\, |u'-u|^2   
\eeqq
\beqq\label{hellinger-13}
E_u\left(  \left[ \wt L^{u' / u} \right]^{1/2}\right)  \;\;\le\;\;   \exp{ \left\{ \,-\; k_\vth\; |u'-u| \right\} }   
\eeqq
with constants $c_\vth$, $k_\vth$ which do not depend on $u$, $u'$. 

\vskip0.5cm
{\bf Proof: }  Again we suppress the scaling factor. We use (\ref{equivariance-1}) to write for $u'\neq u$ 
$$
\left[ \wt L^{u' / u} \right]^{1/2}   \;=\;  \exp{ \left\{ \;\frac12\; \wt{\wt W}(u'-u) \;-\; \frac18\, |u'-u| \right\} }\; e^{-\; \frac18\, |u'-u|} \;. 
$$
By 5.3, the expectation of the first term on the right hand side under $\wt P_u$ is~$1$: this proves (\ref{hellinger-13}) and 
$$ %\beqq\label{hellinger-11-bis}
E_u\left( \left[ 1 -  \left( \wt L^{u' / u} \right)^{1/2} \right]^2 \right)      \;\;=\;\; 2\left( 1 -  e^{-\; \frac18\, |u'-u|} \right)  
$$ %\eeqq
and  (\ref{hellinger-11}). Similiar calculations considering $\left[ \wt L^{u' / u} \right]^{j/4}$ for $j=1,2,3$ give 
$$
E_u\left( \left[ 1 - \left( \wt L^{u' / u} \right)^{1/4} \right]^4 \right) 
\;\;\ = \;\;  2 \;+\; 6\, e^{-\; \frac18\, |u'-u|} \;-\; 8\, e^{-\; \frac{3}{32}\, |u'-u|}  
$$
which behaves as $cst \cdot  |u'-u|^2$ as $u'\to u$. This proves (\ref{hellinger-12}). \halmos \\

Since the left hand side of (\ref{hellinger-11}) is twice the squared Hellinger distance, the parametrization is H\"older continuous with index~$\frac12$ 
\beqq\label{hellingerlimitexperiment}
\lim_{u'\to u}\; \frac{1}{\sqrt{|u'-u|}}\; H\left( \wt P_{u'}  , \wt P_u \right) \;\;=\;\; \sqrt{\frac18 \, J_\vth\;}  
\eeqq
at every point $u$ of the limit experiment $\wt\cale = \left\{ \wt P_u : u\in\bbr \right\}$ of (\ref{limitexperiment})+(\ref{limitlikelihoodratio}) for local models at $\vth$.  \\
%%%%%%%
% diese zwei lemmata veraendert 05.11.09
%%%%%%%%

In order to prove theorem 1.4, we shall follow Ibragimov and Khasminskii [IH 81, theorems 5.1+5.2 in section 1 and 19--21 in appendix A1.4]. \\

{\bf 5.5 Proof of theorem 1.4: } Fix $\vth\in\Theta$; in analogy to [IH 81], we use notations 
$$
Z_{n,\vth}(u)  \;:=\;  L_{nT}^{ (\vth+\frac{u}{n}) / \vth } \;,\; u\in\Theta_{\vth,n} 
 \quad,\quad 
\wt Z (u) \;:=\;  \wt L ^{u / 0} \;,\; u\in\bbr \;. 
$$
%and extend $\Theta_{\vth,n}$ to $\bbr$ by using the default definition indicated in the beginning of this section whenever necessary. For local models $\cale_n^{(\vth)}$ at $\vth\in\Theta$,  
By (\ref{hellinger-2}), there is some $q'\in\bbn$ and some constant $C$ (both not depending on $\vth$ or $n$) such that 
\beqq\label{IHconditon-1}
E_\vth \left( \left[\,  Z^{1/4}_{n,\vth}(u_1)  \;-\;   Z^{1/4}_{n,\vth}(u_2)\,\right]^4 \right)  \;\;\le\;\;  C\, (1+K^{q'})\, |u_1-u_2|^2
\quad\mbox{if  $\;|u_1|,|u_2|\le K$ }
\eeqq
holds for all $K\in\bbn$, $\vth\in\Theta$, $n\in\bbn$.  By theorem 1.1 a), we have convergence as $\nto$  
\beqq\label{IHconditon-2}
\left(  Z^{1/4}_{n,\vth}(u) \right)_{ u\in\Theta_{\vth,n} }   \;\;\lra\;\;    \left(  \wt Z^{1/4}(u) \right)_{u\in\bbr} 
\eeqq
in the sense of finite dimensional distributions. (\ref{hellinger-3}) gives some $k$ (not depending on $\vth$ or $n$) such that 
\beqq\label{IHconditon-3}
E_\vth \left( \,  Z^{1/2}_{n,\vth}(u) \right)  \;\;\le\;\;  e^{\,-\, k\, |u| }
\quad,\quad u\in\Theta_{\vth,n}
\eeqq 
for all $\vth\in\Theta$ and $n\in\bbn$. In  the limit experiment $\wt\cale$, we have the corresponding assertions  
\beqq\label{IHconditon-11}
E_{\wt P_0} \left( \left[\,  \wt Z^{1/4}(u_1)  \;-\;   \wt Z^{1/4}(u_2)\,\right]^4 \right)  \;\;\le\;\;  C\, |u_1-u_2|^2 \quad,\quad u_1,u_2\in\bbr
\eeqq
\beqq\label{IHconditon-12}
E_{\wt P_0} \left( \,  \wt Z^{1/2}(u) \right)  \;\;\le\;\;  e^{\,-\, k\, |u| } \quad,\quad u\in\bbr  
\eeqq
from (\ref{hellinger-12})+(\ref{hellinger-13}) in lemma 5.4. 

1)~ Fix $K<\infty$ arbitrarily large. In virtue of [IH 81, p.\ 378], the assertions (\ref{IHconditon-1})+(\ref{IHconditon-2}) imply weak convergence in $C([{-}K,K])$ of $\,Z^{1/4}_{n,\vth}$ under $Q^\vth$ to $\,\wt Z^{1/4}$ under $\wt P_0$ as $\nto$, and thus weak convergence in $C([{-}K,K])$ of $\,Z_{n,\vth}$  under $Q^\vth$ to $\,\wt Z$ under $\wt P_0$ as $\nto$. This proves part a) of the theorem. 

2)~ Exactly on the lines of  the arguments given in [IH 81,  p.\ 372 and p.\ 42--44], we deduce from assertions (\ref{IHconditon-1})+(\ref{IHconditon-3}) above the following:  
there is  some $q\in\bbn$ and --~with respect to any $K_0>0$  fixed~--  suitable constants  $b_3, b_2 >0$ such that 
\beqq\label{IHconditon-4}
\left\{ \begin{array}{l}
Q^\vth\left( \sup\limits_{K+r\le |u| < K+r+1}\, Z_{n,\vth}(u) \;\ge \vep \; \right)  \;\;\le\;\; \frac{1}{\sqrt{\vep}}\, b_3\, (K{+}r{+}1)^q\, e^{\,-\; b_2 (K+r)}  
\\ \mbox{for any choice of}\quad  \vep>0 \;,\; r\in\bbn_0 \;,\;  K\ge K_0 \;.   
\end{array} \right. 
\eeqq
The constants  $b_3$, $b_2$, $q$ in (\ref{IHconditon-4}) do not depend on $n$ or $\vth$ (since they come from the right hand sides in (\ref{IHconditon-1})+(\ref{IHconditon-3}) and from choice of $K_0$). Similiarly, at the level of the limit experiment $\wt\cale$, we deduce in the same way from  (\ref{IHconditon-11})+(\ref{IHconditon-12}) above  
\beqq\label{IHconditon-4-limitexperiment}
\left\{ \begin{array}{l}
\wt P_0 \left( \sup\limits_{K+r\le |u| < K+r+1}\, \wt Z (u) \;\ge \vep \; \right)  \;\;\le\;\; \frac{1}{\sqrt{\vep}}\, b_3\, (K{+}r{+}1)^q\, e^{\,-\; b_2 (K+r)}  
\\ \mbox{for any choice of}\quad  \vep>0 \;,\; r\in\bbn_0 \;,\;  K\ge K_0  \;. 
\end{array} \right. 
\eeqq

3)~ From (\ref{IHconditon-4}) with $\vep=1$, %we draw a first conclusion: 
we have for all $\vth\in\Theta$ and $n\ge 1$  
\beqq\label{IHconditon-5}
Q^\vth\left( \sup_{ |u| > K }\, Z_{n,\vth}(u) \;\ge 1\; \right)  \;\;\le\;\; b_1\, e^{\,-\; b_2\, K} 
 \;,\quad  K\ge K_0
\eeqq
with constants which do not depend on $\vth$, $n$ or $K\ge K_0$. This is seen similiar to [IH 81,  p.\ 43] after summation over $r\in\bbn_0$ in (\ref{IHconditon-4})  with $\vep=1$. From (\ref{IHconditon-4}) where $\vep$ is replaced by $\vep_r := \frac{1}{(K+r+1)^p}$ which depends on $r$, we obtain assertion (\ref{IHconditon-6}) of theorem 1.4 since 
$$
Q^\vth\left( \sup_{K+r\le |u| < K+r+1}\, |u|^p\, Z_{n,\vth}(u) \;\ge 1\; \right)  \;\;\le\;\; 
Q^\vth\left( \sup_{K+r\le |u| < K+r+1}\, Z_{n,\vth}(u) \;\ge \; \vep_r \; \right)  \;. 
$$
At the level of the limit experiment $\wt\cale$, we start  from (\ref{IHconditon-4-limitexperiment}) and obtain  in  the same way 
\beqq\label{IHconditon-15}
\wt P_0 \left( \sup_{ |u| > K }\, \wt Z(u) \;\ge 1\; \right)  \;\;\le\;\; b_1\, e^{\,-\; b_2\, K} \;,\quad  K\ge K_0
\eeqq
together with  assertion (\ref{IHconditon-16}) of theorem 1.4.  Part b) of theorem 1.4 is proved. 

4)~ Next we consider %--~following [IH 81, p.~45]~-- 
denominators in (\ref{IHconditon-86}) or (\ref{IHconditon-96}) and show that there is some $D>0$ such that 
$$
Q^\vth \left( \int_{\Theta_{\vth,n}} Z_{n,\vth}(u') du'  \,<\,  \frac{\delta}{2}\, \right) \;<\; D\, \sqrt{\delta} 
\quad,\quad 
\wt P_0 \left( \int_{\bbr} \wt Z(u') du'  \,<\, \frac{\delta}{2}\,  \right) \;<\; D\, \sqrt{\delta} 
$$
holds for all $0<\delta<1$ and for all $n , \vth$. As in [IH 81, p.\ 45-46], this comes from 
\beao
&&Q^\vth \left( \int_{ -\frac{\delta}{2} }^{ \frac{\delta}{2} } Z_{n,\vth}(u') du'  \,<\, \frac{\delta}{2}\, \right)  
\;=\; Q^\vth \left( \int_{ -\frac{\delta}{2} }^{ \frac{\delta}{2} } [Z_{n,\vth}(u')-1] du'  \,<\, -\frac{\delta}{2}\, \right)   \\
&&\le\; Q^\vth \left( \int_{ -\frac{\delta}{2} }^{ \frac{\delta}{2} } |Z_{n,\vth}(u')-1| du'  \,>\, \frac{\delta}{2}\, \right)  
\;\;\le\;\;  2\cdot \frac{1}{\delta}  \int_{ -\frac{\delta}{2} }^{ \frac{\delta}{2} } E_\vth\left( |Z_{n,\vth}(u')-1| \right) du'   
\eeao
where the last integrand is a total variation distance, thus smaller than Hellinger distance: hence (\ref{hellinger-1}) for values $|u'|\le\frac12$ of the local parameter gives $E_\vth\left( |Z_{n,\vth}(u')-1| \right) \le cst \sqrt{|u'|}$ and proves the first of the two assertions. The second one follows in the same way from (\ref{hellinger-11}). 

5)~ Next we show on the lines of  [IH 81, p.~47] that for $K_0$ fixed and  all $K\ge K_0$, $r\in\bbn_0$
\beqq\label{lastassertion-1}
E_\vth \left(\; \int_{ \Theta_{\vth,n} \cap \{K+r\le |u|<K+r+1\} }  \frac{ Z_{n,\vth}(u) }{ \int_{\Theta_{\vth,n}} Z_{n,\vth}(u') du'  } \; du\;\right) 
\;\;\le\;\;  cst\cdot e^{\,-\; \frac14\, b_2 (K+r)}
\eeqq
\beqq\label{lastassertion-2}
E_{\wt P_0} \left(\; \int_{  \{K+r\le |u|<K+r+1\} }  \frac{ \wt Z(u) }{ \int_{\bbr} \wt Z(u') du'  } \; du\;\right) 
\;\;\le\;\;  cst\cdot e^{\,-\; \frac14\, b_2 (K+r)}
\eeqq
with respect to the same spherical sections as  in step 2), and with constants which do not depend on $n$ or $\vth$. To see  this, write 
$$
I_r \;:=\; \int_{ \Theta_{\vth,n} \cap \{K+r\le |u|<K+r+1\} } Z_{n,\vth}(u)\, du    \quad,\quad 
H_r \;:=\;  \int_{ \Theta_{\vth,n} \cap \{K+r\le |u|<K+r+1\} }  \frac{ Z_{n,\vth}(u) }{ \int_{\Theta_{\vth,n}} Z_{n,\vth}(u') du'  } \; du \;. 
$$
Then $\{ I_r > 2\vep \}$ is a subset of $\;\{ \sup\limits_{K+r\le |u|<K+r+1} Z_{n,\vth}(u) \ge  \vep \}\,$. Applying (\ref{IHconditon-4}) with  $b_2>0$ as there and  with  $\vep = \vep_r = \frac12 e^{-b_2 (K+r)}$ depending on $r$, we obtain 
$$
Q^\vth \left(\; I_r > e^{-b_2 (K+r)} \;\right)  \;\;\le\;\;  \wt b_3\, (K{+}r{+}1)^q\, e^{\,-\; \frac12\, b_2\, (K+r)}  
\quad,\quad  r\in\bbn_0 \;,\; K\ge K_0 \;. 
$$
From this, using $H_r\le 1$ and step 4), we have for any $0<\delta<1$ 
\beao
E_\vth \left(\, H_r \,\right)  &\le&  Q^\vth \left( \int_{\Theta_{\vth,n}} Z_{n,\vth}(u') du'  \,<\,  \frac{\delta}{2}\, \right) \;+\; Q^\vth \left(\; I_r > e^{-b_2 (K+r)} \;\right)  \;+\; 2\, \frac{1}{\delta}\, e^{\,-\; b_2\, (K+r)}  \\
&\le& cst \cdot \left( \sqrt{\delta}  \;+\;  (K{+}r{+}1)^q\, e^{\,-\; \frac12\, b_2\, (K+r)} \;+\;  \frac{1}{\delta}\, e^{\,-\; b_2\, (K+r)}  \right) \;. 
\eeao
From this, the choice $\delta = \delta_r = e^{\,-\; \frac12\, b_2 (K+r)}$ yields assertion (\ref{lastassertion-1}). The proof of (\ref{lastassertion-2}) is similiar. 

6)~ Finally, we replace the integrals in (\ref{IHconditon-86}) or (\ref{IHconditon-96}) 
$$
E_\vth \left(\; \int_{ \Theta_{\vth,n} \cap \{|u|>K\} }  |u|^p\; \frac{ Z_{n,\vth}(u) }{ \int_{ \Theta_{\vth,n} } Z_{n,\vth}(u') du'  } \; du\;\right) 
\quad\mbox{or}\quad E_{\wt P_0} \left(\; \int_{ \{|u|>K\} }  |u|^p\;  \frac{ \wt Z(u) }{ \int_{\bbr} \wt Z(u') du'  } \; du\;\right) 
$$
for arbitrary $p\in\bbn$ and $K\ge K_0$ by sums 
$$
\sum_{r\in\bbn_0} (K+r+1)^p\;  E_\vth \left(\; \int_{ \Theta_{\vth,n} \cap \{K+r\le |u|<K+r+1\} }  \frac{ Z_{n,\vth}(u) }{ \int_{\Theta_{\vth,n}} Z_{n,\vth}(u') du'  } \; du\;\right) 
$$
or  
$$
\sum_{r\in\bbn_0} (K+r+1)^p\;  E_{\wt P_0} \left(\; \int_{ \{K+r\le |u|<K+r+1\} }  \frac{ \wt Z(u) }{ \int_{\bbr} \wt Z(u') du'  } \; du\;\right)   
$$
and apply (\ref{lastassertion-1}) and (\ref{lastassertion-2}). This finishes the proof of part c) of theorem 1.4.\halmos\\

Bounds (\ref{IHconditon-4})+(\ref{IHconditon-4-limitexperiment}) above control maxima of the likelihood  over spherical regions in the parameter space with center $\vth$. Together with weak convergence of likelihoods, they  imply according to [IH 81] convergence of both maximum likelihood and Bayes estimators together with moments of arbitrary order.  See also the approach of Strasser [S~85, Theorem 65.5 and Corollary 67.6] 
to convergence of Bayes estimators, based on a uniform integrability condition which is satisfied through (\ref{IHconditon-86}) or (\ref{IHconditon-96}). A short exposition of the key to convergence of  maximum likelihood (\ref{MLEstage-n})+(\ref{MLElimitexperiment}) and  Bayes estimators (\ref{BE})+(\ref{BElimitexperiment}) is given in [K~08].\\

%Weak convergence 1.4 a) of likelihoods combined with bounds 1.4 b)+c) leads to convergence of  maximum likelihood estimators (\ref{MLEstage-n})+(\ref{MLElimitexperiment}) and of Bayes estimators (\ref{BE})+(\ref{BElimitexperiment}).  
%The core of this argument, see an exposition in [K 08], is well known.\\

{\bf 5.6 Proof of theorem 1.5 a): }  1)~ From Ibragimov and Khasminskii [IH 81, lemma 2.5 on p.\ 335--336], the following is known for  the MLE in the limit experiment $\wt\cale$: the law $\call (\wh u | \wt P _0 )$ has no point masses, and the argmax in (\ref{MLElimitexperiment}) is unique almost surely. Symmetry in law of two-sided Brownian motion around~$0$ in (\ref{limitlikelihoodratio}) implies that $\call (\wh u | \wt P _0 )$ is symmetric around~$0$. 
Due to (\ref{IHconditon-6})+(\ref{IHconditon-16}), we can choose  $K$ large enough to make 
\beao
Q^\vth\left( |n\, ( \wh\vth_{nT} - \vth ) | > K  \right)  
\;\le\; 
Q^\vth\left(  \sup_{|u|>K}\, L_{nT}^{(\vth + \frac{u}{n}) / \vth }  \ge  \sup_{|u|\le K}\, L_{nT}^{(\vth + \frac{u}{n}) / \vth } \right)  
\;\le\; 
Q^\vth\left(  \sup_{|u|>K}\, L_{nT}^{(\vth + \frac{u}{n}) / \vth }  \ge 1 \right) 
\eeao
arbitrarily small, uniformly in $n$, together with  
$$
\wt P_0\left( |\wh u|>K \right)  \;\le\;  \wt P_0\left(  \sup_{|u|>K}\wt L^{u / 0}  \ge  \sup_{|u|\le K}\, \wt L^{u / 0}  \right) \;\le\; \wt P_0\left(  \sup_{|u|>K}\wt L^{u / 0}  \ge 1 \right) \;. 
$$ 
This allows to approximate for fixed $x$  
$$ %\beqq\label{fact-21}
Q^\vth\left(  n\, ( \wh\vth_{nT} - \vth ) \le x  \right)  
\;=\;  Q^\vth \left(  \;\sup\limits_{u\le x}\,  L_{nT}^{(\vth + \frac{u}{n}) / \vth }    \;\ge\;   \sup\limits_{u\ge x}\,  L_{nT}^{(\vth + \frac{u}{n}) / \vth }  \right)  
$$ %\eeqq
uniformly in $n$ by 
$$ %\beqq\label{fact-23}
Q^\vth \left(  \;\sup\limits_{u\in[-K,x]}\,  L_{nT}^{(\vth + \frac{u}{n}) / \vth }   \;\ge\;   \sup\limits_{u\in[x,K]}\, L_{nT}^{(\vth + \frac{u}{n}) / \vth } \right)   \;, 
$$ %\eeqq
and 
$$ %\beqq\label{fact-22}
\wt P _0 \left(\; \wh u \le x \;\right) \;=\; \wt P _0 \left(  \;\sup\limits_{u\le x}\,  \wt L^{ u / 0 }    \;\ge\;   \sup\limits_{u\ge x}\,  \wt L^{ u / 0 } \right)   
$$ %\eeqq
by 
$$ %\beqq\label{fact-24}
\wt P _0 \left(  \;\sup\limits_{u\in[-K,x]}\,  \wt L^{ u / 0 }    \;\ge\;   \sup\limits_{u\in[x,K]}\,  \wt L^{ u / 0 } \right) \;.  
$$ %\eeqq
The continuous mapping theorem and theorem 1.4~a) give weak convergence 
$$
\call\left(  \;\sup\limits_{u\in[-K,x]}\,  L_{nT}^{(\vth + \frac{u}{n}) / \vth }   ,\;   \sup\limits_{u\in[x,K]}\, L_{nT}^{(\vth + \frac{u}{n}) / \vth } \;\mid\; Q^\vth \right)  
\;\;\lra\;\;
\call\left(  \;\sup\limits_{u\in[-K,x]}\,  \wt L^{ u / 0 }   ,\;   \sup\limits_{u\in[x,K]}\,  \wt L^{ u / 0 } \;\mid\; \wt P_0\right)  
$$
in $\bbr^2$ as $\nto$, where  $0$ is a continuity point for the law of the difference 
$$
\call\left(  \;\sup\limits_{u\in[-K,x]}\,  \wt L^{ u / 0 } \;-   \sup\limits_{u\in[x,K]}\,  \wt L^{ u / 0 }  \mid \,\wt P_0\, \right)  
$$
again by [IH 81, lemma 2.5 on p.\ 335--336]. This is weak convergence of MLE.

2)~ Writing rescaled estimation errors at $\vth$ for Bayes estimators (\ref{BE}) as  
\beqq\label{rescBEerrors}
n\;  \frac{\int_{\Theta} (\zeta-\vth) \,  L_{nT}^{\zeta / \zeta_0 }\, d\zeta}{ \int_{\Theta}  L_{nT}^{\zeta / \zeta_0 }\, d\zeta }  
\;\;=\;\;     n\;  \frac{\int_{\Theta} (\zeta-\vth)\,  L_{nT}^{\zeta / \vth }\, d\zeta}{ \int_{\Theta}  L_{nT}^{\zeta / \vth }\, d\zeta }  
\;\;=\;\;    \frac{\int_{\Theta_{\vth,n}} u\;  L_{nT}^{\vth+\frac{u}{n} / \vth }\, du}{ \int_{\Theta_{\vth,n}}  L_{nT}^{\vth+\frac{u}{n} / \vth }\, du }   
\eeqq
we  have   
$$ %\beqq\label{fact-31}
Q^\vth\left(  n\, ( \vth^*_{nT} - \vth ) \le x  \right)  
\;\;=\;\;    Q^\vth \left(\;    \frac{\int_{\Theta_{\vth,n}} u\;  L_{nT}^{\vth+\frac{u}{n} / \vth }\, du}{ \int_{\Theta_{\vth,n}}  L_{nT}^{\vth+\frac{u}{n} / \vth }\, du } \le x  \;\right)    
$$ %\eeqq
together with 
$$ %\beqq\label{fact-32}
\wt P_0 \left(\; u^* \le x  \; \right)
\;\;=\;\; 
\wt P_0 \left(\;   \frac{ \int_{-\infty}^\infty u\,  \wt L^{ u / 0 }\, du }{ \int_{-\infty}^\infty  \wt L^{ u / 0 }\, du }  \le x   \; \right)   
$$ %\eeqq
for (\ref{BElimitexperiment}) in the limit experiment. From (\ref{IHconditon-6}) with $p=3$ we get for all $K\ge K_0$  
\beao
Q^\vth \left( \int_{ \Theta_{\vth,n} \cap \{|u|>K\}} |u|\, L_{nT}^{ (\vth+\frac{u}{n}) \,/\,  \vth } du\;  > \frac{1}{K} \right) 
\;\le\; 
Q^\vth \left(  L_{nT}^{ (\vth+\frac{u}{n}) \,/\,  \vth }  \ge  |u|^{-3} \;\;\mbox{for some $|u|\ge K$}  \right) 
\;\le\;  b_1\, e^{- b_2 K }  
\eeao
independently of $n$; similiarly, in the limit experiment, (\ref{IHconditon-16}) with $p=3$ gives 
$$
\wt P_0 \left(\,  \int_{\{|u|>K\}} |u|\, \wt L^{ u \,/\,  0 }\, du   \;>\; \frac{1}{K}\,  \right) \;\; \le \;\;  b_1\, e^{- b_2 K } 
$$
for all $K\ge K_0$. Chosing $K$ large enough, we approximate 
\beqq\label{term-n-original}
\call \left( \int_{\Theta_{\vth,n}} u\;  L_{nT}^{\vth+\frac{u}{n} / \vth }\, du \;,\; \int_{\Theta_{\vth,n}} L_{nT}^{\vth+\frac{u}{n} / \vth }\, du \;\mid\; Q^\vth \right)  
\eeqq
by 
\beqq\label{term-n}
\call \left( \int_{[-K,K]} u\;  L_{nT}^{\vth+\frac{u}{n} / \vth }\, du \;,\; \int_{[-K,K]} L_{nT}^{\vth+\frac{u}{n} / \vth }\, du \;\mid\; Q^\vth \right)   
\eeqq
such that the accuracy of this approximation does not depend on $n$,  and similiarly %in the limit experiment 
\beqq\label{term-unendlich-original}
\call \left(  \int_{-\infty}^\infty u\,  \wt L^{ u / 0 }\, du \;,\;  \int_{-\infty}^\infty \wt L^{ u / 0 }\, du \;\mid\; \wt P_0 \right)  
\eeqq
by 
\beqq\label{term-unendlich}
\call \left(  \int_{[-K,K]}  u\,  \wt L^{ u / 0 }\, du \;,\;  \int_{[-K,K]}  \wt L^{ u / 0 }\, du \;\mid\; \wt P_0 \right) \;. 
\eeqq
From 1.4~a) and the continuous mapping theorem, we have weak convergence in $\bbr^2$  as $\nto$ of (\ref{term-n})  to (\ref{term-unendlich}), for arbitrary $K\ge K_0$. This yields weak convergence of (\ref{term-n-original}) to  (\ref{term-unendlich-original}). The second component in (\ref{term-unendlich-original}) being strictly positive, we have weak convergence of BE. \halmos\\

{\bf 5.7 Proof of theorem 1.5 b): } For arbitrary  $H_0>0$ fixed, there are constants $\wt b_1 , \wt b_2$  such that 
\beqq\label{IHconditon-21}
\sup_n\; Q^\vth \left(\, \left| n\,(\vth^*_{nT}- \vth)\right| > H \,\right)  \;\;\le\;\;  \wt b_1\, e^{ - \wt b_2\, H} 
\;,\quad H>H_0   
\eeqq
holds: from (\ref{rescBEerrors}) and the trivial 
$$
\int_{ \Theta_{n,\vth} \cap \{|u|\le K\} } |u|\, \frac{ L_{nT}^{ (\vth+\frac{u}{n}) / \vth } }{ \int_{ \Theta_{n,\vth} }  L_{nT}^{ (\vth+\frac{u'}{n}) / \vth } du' }\, du  \;\;\le\;\;  K
$$ 
we get 
\beao
Q^\vth \left(  \int_{ \Theta_{n,\vth} } |u|\, \frac{ L_{nT}^{ (\vth+\frac{u}{n}) / \vth } }{ \int_{ \Theta_{n,\vth} }  L_{nT}^{ (\vth+\frac{u'}{n}) / \vth } du' }\, du  \,>\, 2K  \right)  
\;\;\le\;\; 
\frac1K\; E_\vth \left( \int_{ \Theta_{n,\vth} \cap \{|u|>K\} } |u|\, \frac{ L_{nT}^{ (\vth+\frac{u}{n}) / \vth } }{ \int_{ \Theta_{n,\vth} }  L_{nT}^{ (\vth+\frac{u'}{n}) / \vth } du' }\, du \right)   
\eeao
and apply  (\ref {IHconditon-86}). For MLE, as in the beginning of step 1) of 5.6, we have 
\beqq\label{MLtightness}
Q^\vth \left(\, \left| n\,(\wh \vth_{nT}- \vth)\right| > H \,\right)  \;\le\;  Q^\vth\left( \sup_{ |u| > H }\, Z_{n,\vth}(u) \;\ge 1\; \right)  \;\;\le\;\; b_1\, e^{\,-\; b_2\, H}  \;,\quad  H\ge H_0   
\eeqq
by (\ref{IHconditon-6}) or (\ref{IHconditon-5}). By weak convergence of rescaled estimation errors established in 5.6, bounds (\ref{IHconditon-21})+(\ref{MLtightness}) carry over to the limit experiment and give  
\beqq\label{IHconditon-21-bis}
\wt P_0 \left( |u^*| > H \right) \;\le\;  \wt b_1\, e^{-\wt b_2\, H} \quad,\quad 
\wt P_0 \left( |\wh u| > H \right) \;\le\;  b_1\, e^{-b_2\, H} \quad,\quad H\ge H_0 \;. 
\eeqq
Now convergence of moments of order $p\in\bbn$ for MLE and BE  
$$
\int_0^\infty  y^{p-1}\; Q^\vth \left(\, \left| n\,(\wh \vth_{nT}- \vth)\right| > y \,\right)\; dy   \quad,\quad  
\int_0^\infty  y^{p-1}\; Q^\vth \left(\, \left| n\,(\vth^*_{nT}- \vth)\right| > y \,\right)\; dy
$$
as $\nto$ to 
$$
\int_0^\infty y^{p-1}\; \wt P_0 \left(\, \left| \wh u \right| > y \,\right)\; dy \quad,\quad  
\int_0^\infty y^{p-1}\; \wt P_0 \left(\, \left| u^* \right| > y \,\right)\; dy
$$
is a consequence of weak convergence of recaled estimation errors in combination with dominated convergence thanks to (\ref{IHconditon-21})+(\ref{MLtightness})+(\ref{IHconditon-21-bis}).\halmos\\

%Now we discuss equivariance properties of the BE sequence and prove proposition 1.7. \\

{\bf 5.8 Proof of theorem 1.7: } 1)~ By lemma 5.3, $\;\left\{ \wt P_{u+h} : h\in\bbr \right\}$ being statistically the same experiment as $\wt\cale = \left\{ \wt P_u : u\in\bbr \right\}$, laws of rescaled estimation errors  
$$
\call \left( u^*-u \mid \wt P_u \right) \;\;=\;\; 
\call \left(\; \frac{ \int_{-\infty}^\infty u'\,  \wt L^{ u' / 0 }\, du' }{ \int_{-\infty}^\infty  \wt L^{ u' / 0 }\, du' } - u \;\mid \wt P_u \right) 
\;\;=\;\; 
\call \left(\; \frac{ \int_{-\infty}^\infty h\,  \wt L^{ (u+h) / u }\, dh }{ \int_{-\infty}^\infty  \wt L^{ (u+h) / u }\, dh }  \;\mid \wt P_u \right)     
$$
and 
$$
\call \left( \wh u -u \mid \wt P_u \right) \;\;=\;\; 
\call \left(\; \mathop{\rm argmax}\limits_{h\in\bbr} \wt L^{(u+h)/u} \mid \wt P_u \right)  
$$
do not depend on $u\in\bbr$. Hence both BE (\ref{BElimitexperiment}) and MLE (\ref{MLElimitexperiment}) are equivariant estimators for the parameter in the limit experiment $\wt\cale$. Thus 1.7~a) holds.

2)~  In order to prove 1.7~b), we shall use 'LeCam's Third Lemma' for contiguous alternatives (see [LY 90, pp.\ 22--23], or [H 08, 3.6+3.16]) in combination with the above equivariance property of the BE in the limit experiment. We have seen in the proof 5.6 that pairs (\ref{term-n-original}) %under $Q^\vth$ 
$$
\left( \int_{\Theta_{\vth,n}} u\;  L_{nT}^{\vth+\frac{u}{n} / \vth }\, du \;,\; \int_{\Theta_{\vth,n}} L_{nT}^{\vth+\frac{u}{n} / \vth }\, du\right)  \quad\mbox{under $Q^\vth$}
$$
converge weakly in $\bbr^2$ as $\nto$ to  the pair (\ref{term-unendlich-original}) %under $\wt P_0$ 
$$
\left(  \int_{-\infty}^\infty u\,  \wt L^{ u / 0 }\, du \;,\;  \int_{-\infty}^\infty \wt L^{ u / 0 }\, du \right) 
\quad\mbox{under $\wt P_0$} \;. 
$$
Obviously, for $u_0\in\bbr$ fixed, the same argument also yields joint convergence of triplets 
\beqq\label{term-n-original-2}
\left(\; L_{nT}^{\vth+\frac{u_0}{n} / \vth } \;,\; \int_{\Theta_{\vth,n}} u\;  L_{nT}^{\vth+\frac{u}{n} / \vth }\, du \;,\; \int_{\Theta_{\vth,n}} L_{nT}^{\vth+\frac{u}{n} / \vth }\, du\right)  \quad\mbox{under $Q^\vth$}
\eeqq
weakly in $\bbr^3$ as $\nto$ to 
$$
\left(\; \wt L^{ u_0 / 0 } \;,\; \int_{-\infty}^\infty u\,  \wt L^{ u / 0 }\, du \;,\;  \int_{-\infty}^\infty \wt L^{ u / 0 }\, du \right)  \quad\mbox{under $\wt P_0$} \;. 
$$
For any convergent sequence $(u_n)_n$ tending to the limit $u_0$, this convergence  remains valid if we  place $L_{nT}^{\vth+\frac{u_n}{n} / \vth }$ instead of $L_{nT}^{\vth+\frac{u_0}{n} / \vth }$ into the first component of (\ref{term-n-original-2}). From joint convergence with the sequence of likelihood ratios, LeCam's Third Lemma deduces weak convergence under the corresponding contiguous alternatives:  thus 
$$
\left(\; L_{nT}^{\vth+\frac{u_n}{n} / \vth } \;,\; \int_{\Theta_{\vth,n}} u\;  L_{nT}^{\vth+\frac{u}{n} / \vth }\, du \;,\; \int_{\Theta_{\vth,n}} L_{nT}^{\vth+\frac{u}{n} / \vth }\, du\right)  \quad\mbox{under $Q^{\vth+\frac{u_n}{n}}$}
$$
converges weakly in $\bbr^3$ as $\nto$ to 
$$
\left(\; \wt L^{ u_0 / 0 } \;,\; \int_{-\infty}^\infty u\,  \wt L^{ u / 0 }\, du \;,\;  \int_{-\infty}^\infty \wt L^{ u / 0 }\, du \right)  \quad\mbox{under $\wt P_{u_0}$} \;. 
$$
Using the continuous mapping theorem we obtain weak convergence of ratios 
$$
\frac{ \int_{\Theta_{\vth,n}} (u-u_n)\;  L_{nT}^{\vth+\frac{u}{n} / 0 }\, du  }{ \int_{\Theta_{\vth,n}} L_{nT}^{\vth+\frac{u}{n} / 0 }\, du } \quad\mbox{under $Q^{\vth+\frac{u_n}{n}}$}
$$
as $\nto$ to 
$$
\frac{ \int_{-\infty}^\infty (u-u_0)\,  \wt L^{ u / 0 }\, du }{  \int_{-\infty}^\infty \wt L^{ u / 0 }\, du } \quad\mbox{under $\wt P_{u_0}$} \;. 
$$
Thus we have proved weak  convergence as $\nto$ of rescaled BE errors 
\beqq\label{rescBEerrors-contiguous}
\call\left(\, n(\vth^*_{nT}-(\vth+\frac{u_n}{n})) \mid Q^{\vth+\frac{u_n}{n}} \right)
\;\;\lra\;\;
\call\left( u^* -u_0 \mid \wt P_{u_0}\right)   
\eeqq
under continguous alternatives, for arbitrary convergent sequences $u_n\to u_0$. 
An argument of the same structure works for rescaled MLE errors.

3)~ Consider a loss function $\ell(\cdot)$ which is continuous, subconvex and bounded. Then (\ref{rescBEerrors-contiguous}) gives  
\beqq\label{fact-41}
E_{\vth+\frac{u_n}{n}} \left(\, \ell\left(  n(\vth^*_{nT}-(\vth+\frac{u_n}{n})) \right) \right)
\;\;\lra\;\;
E_{\wt P_0}\left( \ell\left( u^* -u_0 \right)\right)   
\eeqq
for arbitrary convergent sequences $u_n\to u_0$. Via selection of convergent subsequences in compacts $[-C,C]$, for $C$ arbitrarily large, (\ref{fact-41}) shows 
\beqq\label{fact-42}
\lim_{\nto}\; \sup_{|u|\le C}\; \left|\;  E_{\vth+\frac{u}{n}} \left( \ell \left( n(\vth^*_{nT}-(\vth{+}\frac{u}{n})) \right)\right) 
\;-\;  E_{\wt P_u} \left( \ell \left( u^* - u \right)\right)   \;\right| \quad=\quad 0
\eeqq
for loss functions $\ell(\cdot)$ which are continuous, subconvex and bounded. This is part b) of theorem 1.7,  in the special case of bounded loss functions. However, this is sufficient to prove part b) of 1.7 in general. We have exponential decrease in the bounds (\ref{IHconditon-21}) -- where the constants are independent of $\vth$ and $n$ -- and  polynomial bounds for $\ell$, thus contributions 
$$
\sup_{n\ge n_0}\; \sup_{\vth'\in\Theta} \int_{y>K} y^{p-1}\; Q^{\vth'} \left(\, \left| n\,(\vth^*_{nT}- \vth')\right| > y \,\right)\; dy    
$$
can be made arbitrarily small by suitable choice of $K$. At the level of the limit experiments, we use (\ref{IHconditon-21-bis}) and equivariance of $u^*$ proved in step 1). Combined with (\ref{fact-42}), this finishes the proof of 1.7~b). Again this works similiarly for MLE. \halmos\\

We turn to the proof of theorem 1.8.\\

{\bf 5.9 Proof of theorem 1.8: }  Fix $\vth\in\Theta$. Throughout this proof, we consider the particular loss function $\ell(x)=x^2$. Part b) of theorem 1.8 being immediate from 1.7~b), we have to prove part a).  
The local asymptotic minimax bound in a) will be a consequence of convergence of experiments, of a general asymptotic minimax theorem with respect to a fixed loss function $\ell(\cdot)$ given in Strasser [S~85, Corollary 62.6] (see also LeCam [L~72], Millar [M~83, p.\ 91], van der Vaart [V~91, Theorem 6.1]), and the fact that we consider Bayes estimators with respect to squared loss in the limit experiment $\wt\cale$. 

1)~ By [S~85, Corollary 62.6] or [M~83, p.\ 91], we have for $C<\infty$ arbitrarily large but fixed 
\beam\label{meineletztegleichung-1}
&&\liminf\limits_{\nto}\; \inf_{\wt \vth_{nT}}\; \sup_{ |u|\le C }\; E_{\vth+\frac{u}{n}}\left(\; \left[ n\left( \wt \vth_{nT} - (\vth{+}\frac{u}{n}) \right)  \right]^2 \;\right)   \nonumber \\
&& \ge\quad 
\inf_{\wt u}\; \sup_{ |u|\le C}\; E_{\wt P_u}\left(\; \left[ \wt u - u \right]^2 \;\right) 
\quad\ge\quad 
\inf_{\wt u}\; \int_{-C}^C \frac{du}{2C}\, E_{\wt P_u}\left(\; \left[ \wt u - u \right]^2 \;\right) 
\eeam
where the 'inf' on the left hand side is over all estimators based on observation of the process $\xi$ up to time $nT$, and the 'inf' on the right hand side over all transition probabilities from $(\wt\Omega,\wt\cala)$ to $(\bbr,\calb(\bbr))$ in the limit experiment $\wt \cale$ of (\ref{limitlikelihoodratio})+(\ref{limitexperiment}).  

2)~ For $C$ fixed, define 
$$
u^*_C(\om)  \;:=\;  \frac{ \int_{-C}^C u'\, \wt L^{u'/0}(\om)\, du' }{ \int_{-C}^C \wt L^{u'/0}(\om)\, du' } \;\;,\;\; \om\in\wt\Omega 
%\quad\mbox{taking values in $(-C,C)$}
$$
taking values in $(-C,C)$. Write $\calr_C(du')$ for the uniform law on $(-C,C)$. Consider the probability 
$$
\mathbb{P}_C (du',d\om) \;\;:=\;\; \calr_C(du')\, \wt P_{u'}(d\om) \;\;=\;\; \calr_C(du')\, \wt P_0(d\om)\; \wt L^{u'/0}(\om)
$$
on $\left(\, (-C,C){\times}\wt\Omega \;,\, \calb(-C,C){\otimes}\wt\cala \,\right)$. By the Bayes property, or the $L^2\left(\,(-C,C){\times}\wt\Omega \,,\,\mathbb{P}_C \,\right)$-projection property of conditional expectations, we can continue the inequality of step 1) in the form 
\beqq\label{meineletztegleichung-2}
\inf_{\wt u}\; E_{\mathbb{P}_C}\left( \left[ \wt u - u \right]^2 \right)  \quad\ge\quad 
E_{\mathbb{P}_C}\left( \left[ u^*_C - u \right]^2 \right)  \quad\ge\quad  
E_{\mathbb{P}_C}\left( \left[ \mathbb{V}_C \right]^2 \right)  
\eeqq
where we define (all probability laws $\wt P_{u'}$ in the limit experiment $\wt\cale$ being equivalent) 
$$
\mathbb{V}_C (u,\om) \;\;:=\;\; \frac{ \int_{\bbr} 1_{(-C,C)}(u')\, (u'-u)\, \wt L^{u' / 0}(\om)\, du' }{ \int_{\bbr} \wt L^{u' / 0}(\om)\, du' } 
\;\;=\;\; \int_{\bbr} 1_{(-C-u,C-u)}(h)\, h\, \frac{ \wt L^{(u+h) / u}(\om) }{ \int_{\bbr} \wt L^{(u+h) / u}(\om)\, dh } \, dh   
$$
and use the trivial inequality $\left[ u^*_C - u \right]^2  \ge \left[ \mathbb{V}_C \right]^2$ on  $(-C,C){\times}\wt\Omega$. 

3)~ Given inequalities (\ref{meineletztegleichung-1})+(\ref{meineletztegleichung-2}), the proof of part b) of theorem 1.8 will be finished if we prove  
\beqq\label{meineletztegleichung-3}
\lim_{C\to\infty}\,  E_{\mathbb{P}_C}\left( \left[ \mathbb{V}_C \right]^2 \right)  \;\;=\;\;  E_{\wt P_0}\left( [u^*]^2\right) 
\eeqq
for the Bayes estimator $u^*$ with 'uniform prior on the real line'  
$$
u^*-u \;\;=\;\; \int_{\bbr} (u'-u)\, \frac{ \wt L^{u' / 0}(\om) }{ \int_{\bbr} \wt L^{u' / 0}(\om)\, du' } \, du'  \;\;=\;\; 
\int_{\bbr} \;h\; \frac{ \wt L^{(u+h) / u}(\om) }{ \int_{\bbr} \wt L^{(u+h) / u}(\om)\, dh } \, dh   
$$
in the limit experiment $\wt\cale$, cf.\ (\ref{BElimitexperiment}) and 1.7~a). 

4)~ It remains to prove (\ref{meineletztegleichung-3}). Introducing 
\beqq\label{meineletztegleichung-4bis}
\Phi_{C,u}(h) \;:=\; 1_{(-C-u,C-u)}(h) - 1_{(-C,C)}(h) \quad,\quad u\in (-C,C) \;,  
\eeqq
and
$$
\rho^{(1)}_{C,u} \;:=\; \int_{\bbr} \Phi_{C,u}(h)\; h\; \frac{ \wt L^{(u+h) / u} }{ \int_{\bbr} \wt L^{(u+h) / u}\, dh } \; dh 
\quad,\quad 
\rho^{(2)}_{C,u}  \;:=\;  - \int_{\bbr} 1_{\{|h|\ge C\}}\; h\; \frac{ \wt L^{(u+h) / u} }{ \int_{\bbr} \wt L^{(u+h) / u}\, dh } \; dh    
$$
we have a representation  
\beqq\label{meineletztegleichung-5}
\mathbb{V}_C (u,\om)  \;\;=\;\; \left( u^*(\om)-u \right) \;+\; \rho^{(1)}_{C,u}(\om) \;+\; \rho^{(2)}_{C,u}(\om)   
\quad\mbox{on $(-C,C){\times}\wt\Omega$}
\eeqq
where the law under $\wt P_u$ of the third term on the right hand side of (\ref{meineletztegleichung-5}) does not depend on $u\in\bbr$ (cf.\ lemma 5.3 and the first part of the proof 5.8). Since $u^*$ is equivariant in $\wt\cale$ and has finite variance, it is sufficient to show  
\beqq\label{meineletztegleichung-4}
\lim_{C\to\infty}\,  E_{\mathbb{P}_C}\left(\, \left[ \mathbb{V}_C - (u^*-u) \right]^2 \right)  \;\;=\;\;  0  
\eeqq
to establish (\ref{meineletztegleichung-3}). For $\vep>0$ arbitrarily small, choose first $K=K(\vep)$ large enough for 
$$
E_{\wt P_u}\left( \int_{\bbr} 1_{\{|h|> K\}}\; h^2\; \frac{ \wt L^{(u+h) / u} }{ \int_{\bbr} \wt L^{(u+h) / u}\, dh } \; dh  \right) \;\;<\;\; \vep 
$$
using theorem 1.4~c), and then $C=C(\vep)$ large enough for $C > K$ together with    
$$
\frac{K}{C}\; E_{\wt P_u}\left( \int_{\bbr} h^2\; \frac{ \wt L^{(u+h) / u} }{ \int_{\bbr} \wt L^{(u+h) / u}\, dh } \; dh  \right) \;\;<\;\; \vep \;. 
$$
In both last inequalities, the left hand side does not depend on $u$. For $\rho^{(2)}_{C,u}$, Jensen inequality gives 
$$
\int_{-C}^C \frac{du}{2C}\, E_{\wt P_u}\left( [\rho^{(2)}_{C,u} ]^2 \right)  \;\;\le\;\;   \vep  
$$
since $C>K$. Considering $\rho^{(1)}_{C,u}$, we have the following bounds for $u\in (-C,C)$: first, $\Phi_{C,u} \le 1_{\{|\cdot|> K\}}$ as long as $|u|<C-K$, second,  $\Phi_{C,u} \le 1$ in the remaining cases $C-K\le |u|<C$. Again with Jensen, 
$$
\int_{-C}^C \frac{du}{2C}\, E_{\wt P_u}\left( [\rho^{(1)}_{C,u} ]^2 \right)
$$ 
is thus by our choice of $C$ and $K$ smaller than 
\beao
\int_{-C}^{-C+K} \frac{du}{2C}\, E_{\wt P_u}\left(  [\rho^{(1)}_{C,u} ]^2 \right)  \;+\;   \frac{C-K}{C}\, \vep  \;+\;  \int_{C-K}^{C} \frac{du}{2C}\, E_{\wt P_u}\left(  [\rho^{(1)}_{C,u} ]^2 \right)  
\quad\le\quad   2\,\vep \;. 
\eeao
Combining the last two bounds with (\ref{meineletztegleichung-4bis}), this gives 
\beao
&&E_{\mathbb{P}_C}\left( \left[ \mathbb{V}_C - (u^*-u) \right]^2 \right)  
\;\;\le\;\;  2\, \int_{-C}^C \frac{du}{2C}\, E_{\wt P_u}\left( [\rho^{(1)}_{C,u} ]^2 + [\rho^{(2)}_{C,u} ]^2 \right)
\;\;\le\;\;  6\, \vep   \;. 
\eeao
Since $\vep>0$ was arbitrary, this proves (\ref{meineletztegleichung-4}) and thus (\ref{meineletztegleichung-3}). By (\ref{meineletztegleichung-1})+(\ref{meineletztegleichung-2})+(\ref{meineletztegleichung-3}), the proof  of theorem 1.8 is  finished. \halmos \\

%%%%%%%%%%%%%%%%%
%\newpage
\vskip1.0cm
{\Large\bf References} 

%\small
\vskip0.2cm
[ADR 69]\quad
Az\'ema, J., Duflo, M., Revuz, D.: Mesures invariantes des processus de Markov r\'ecurrents. \\
S\'eminaire de Probabilit\'es III, Lecture Notes in Mathematics {\bf 88}, 24--33. Springer 1969. 

\vskip0.2cm
[B 98]\quad
Bass, R.: Diffusions and elliptic operators. Springer 1998. 

\vskip0.2cm
[B 05]\quad
Brandt, C: Partial reconstruction of the trajectories of a discretely observed branching diffusion with immigration and an application to inference. 
PhD thesis, Universit\"at Mainz 2005. 

\vskip0.2cm
[BH 06]\quad
Brodda, K., H\"opfner, R.: A stochastic model and a functional central limit theorem for information processing in large systems of neurons. 
J.\ Math.\ Biol.\ {\bf 52}, 439--457 (2006). 

\vskip0.2cm 
[CK 09]\quad
Chan, N., Kutoyants, Yu.: On parameter estimation of threshold autoregressive models. Preprint 2009. 

\vskip0.2cm 
[D 09]\quad
Dachian, S.: On limiting likelihood ratio processes of some change-point type statistical models. Preprint 2009, arXiv:0907.0440. 

\vskip0.2cm 
[D 85]\quad
Davies, R.: Asymptotic inference when the amount of information is random. 
In: Le Cam, L., Olshen, R. (Eds): Proc. of the Berkeley Symposium  in honour of J. Neyman and J. Kiefer. Vol. II. Wadsworth 1985. 

\vskip0.2cm 
[DP 84]\quad
Deshayes, J., Picard, D.: Lois asymptotiques des tests et estimateurs de rupture dans un mod\`ele statistique classique. 
Annales I.H.P.\ B.\ {\bf 20}, 309--327 (1984). 

\vskip0.2cm 
[DL 05]\quad
Dithlevsen, S., L\'ansk\'y, P.: Estimation of the input parameters in the Ornstein Uhlenbeck neuronal model. 
Phys.\ Rev.\ E {\bf 71}, 011907 (2005). 

\vskip0.2cm 
[DL 06]\quad
Dithlevsen, S., L\'ansk\'y, P.: Estimation of the input parameters in the Feller neuronal model. 
Phys.\ Rev.\ E {\bf 73}, 061910 (2006). 

\vskip0.2cm
[DZ 01]\quad
Dzhaparidze, K., van Zanten, H.: On Bernstein type inequalities for martingales. \\
Stoch.\ Proc.\ Appl.\ {\bf 93}, 109--117 (2001).  

\vskip0.2cm
[G 79]\quad
Golubev, G.: Computation of efficiency of maximum-likelihood estimate when observing a discontinuous signal in white noise  (Russian).  Problems Inform.\ Transmission  {\bf 15}, 61--69  (1979). 
 
\vskip0.2cm 
[H 70]\quad 
H\'ajek, J.: A characterization theorem of limiting distributions for regular estimators.\\  
Zeitschr.\ Wahrscheinlichkeitstheor.\ Verw.\ Geb.\ {\bf 14}, 323--330, 1970. 
 
\vskip0.2cm 
[HS 67]\quad 
H\'ajek, J., Sid\'ak, Z:  Theory of rank tests. Academic Press 1967. 

\vskip0.2cm 
[H 07]\quad 
H\"opfner, R.: On a set of data for the membrane potential in a neuron. \\
Math.\ Biosci.\ {\bf 207}, 275--301 (2007). 

\vskip0.2cm 
[H 08]\quad 
H\"opfner, R.: Asymptotische Statistik. Manuscript 2008.\\ 
{\small\tt http://www.mathematik.uni-mainz.de/$\sim$hoepfner}

\vskip0.2cm
[IH 81]\quad
Ibragimov, I., Has'minskii, R.: Statistical estimation. Springer 1981. 

\vskip0.2cm
[IW 89]\quad
Ikeda, N., Watanabe, S.: Stochastic differential equations and diffusion processes. 2nd ed.\ North-Holland/Kodansha 1989. 

\vskip0.2cm
[JS 87]\quad
Jacod, J., Shiryaev, A.: Limit theorems for stochastic processes. Springer 1987.

\vskip0.2cm 
[J 82]\quad
Jeganathan, P.: On the asymptotic theory of estimation when he limit of log-likelihoods is mixed normal. 
Sankhy{$\ov{\rm a}$}  A {\bf 44}, 173--212 (1982). 

\vskip0.2cm 
[J 95]\quad
Jeganathan, P.: Some aspects of asymptotic theory with applications to time series models. 
Econometric Theory {\bf 11}, 818--887 (1995),  preprint version 1988. 

\vskip0.2cm
[KS 91]\quad
Karatzas, J., Shreve, S.: Brownian motion and stochastic calculus. 2nd ed.\ Springer 1991. 

\vskip0.2cm
[KK 00]\quad
K\"uchler, U., Kutoyants, Y.: Delay estimation for some stationary diffusion-type processes. \\
Scand.\ J.\ Statist.\ {\bf 27}, 405--414 (2000). 

\vskip0.2cm
[K 04]\quad
Kutoyants, Y.: Statistical inference for ergodic diffusion processes. Springer 2004. 

\vskip0.2cm
[K 08]\quad
Kutoyants, Y.: Guest lectures 'Statistical inference for diffusions' given at the University of Mainz, 2008, 
{\small\tt http://www.mathematik.uni-mainz.de/$\sim$hoepfner/Kutoyants.html}

\vskip0.2cm 
[LL 87]\quad
L\'ansk\'y, P., L\'ansk\'a, V.: Diffusion approximation of the neuronal model with synaptic reversal potentials. 
Biol.\ Cybern.\ {\bf 56}, 19--26 (1987).

[LS 99]\quad
L\'ansk\'y, P., Sato, S.: The stochastic diffusion models of nerve membrane depolarization and interspike interval generation. 
J.\ Peripheral Nervous System {\bf 4}, 27--42 (1999). 

\vskip0.2cm 
[L 68]\quad
Le Cam, L.: Th\'eorie asymptotique de la d\'ecision statistique. Montr\'eal 1969. 

\vskip0.2cm 
[L 72]\quad
Le Cam, L.: Limits of experiments. Proc.\ 6th Berkeley Symposium Math.\ Statist.\ Probability, Vol.\ I, 245--261. Univ.\ California Press, 1972. 

\vskip0.2cm 
[LY 90]\quad
Le Cam, L., Yang, G.: Asymptotics in statistics. Some basic concepts. 
Springer 1990.   %(2nd Ed.\ Springer 2002).

\vskip0.2cm
[LS 81]\quad
Liptser, R., Shiryaev, A.: Statistics of random processes, Vols.\ I+II.\\
Springer 1981, 2nd Ed.\ 2001. 

\vskip0.2cm
[MT 93]\quad
Meyn, S., Tweedie, R.: Markov chains and stochastic stability. Springer 1993. 

\vskip0.2cm
[M 83]\quad
Millar, P.: The minimax principle in asymptotic statistical theory. Ecole d'\'et\'e de probabilit\'es de Saint-Flour XI (1981). LNM 976, Springer 1983. 

\vskip0.2cm
[R 75]\quad
Revuz, D.: Markov chains. North Holland 1975.  

\vskip0.2cm
[RY 91]\quad
Revuz, D., Yor, M.: 
Continuous martingales and Brownian motion. Springer 1991. 

\vskip0.2cm
[RS 95]\quad
Rubin, H., Song, K.: Exact computation of the asymptotic efficiency of maximum likelihood estimators of a dicontinuous signal in a Gaussian white noise. Ann.\ Statist.\ {23}, 732--739 (1995).

\vskip0.2cm
[S 85]\quad
Strasser, H.: Mathematical theory of  statistics. de Gruyter 1985. 

\vskip0.2cm
[T 68]\quad
Terent'yev, A.: Probability distribution of a time location of an absolute maximum at the output of a synchronized filter. Radioengineering and Electronics {\bf 13}, 652--657 (1968). 

\vskip0.2cm
[T 89]\quad
Tuckwell, H.: Stochastic processes in the neurosciences. 
CBMS-NSF conference series in applied mathematics, SIAM 1989.

\vskip0.2cm
[V 91]\quad
van der Vaart, A.: An asymptotic represententation theorem. \\
Intern.\ Statist.\ Review {\bf 59}, 97-121 (1991).

%%%%%%%%%%%%%%%%%%%%%%%

%%%%%%%%%%%%%%%%%%%%%
 
\vskip0.5cm
~\hfill {\bf 22.02.2010}

\small
Yury A.\ Kutoyants\\
Laboratoire de Statistique et Processus, Universit\'e du Maine, F--72085 Le Mans Cedex 9\\ 
{\tt kutoyants@univ-lemans.fr}\\
{\tt http://www.univ-lemans.fr/sciences/statist/pages$\_$persos/kuto.html}

Reinhard H\"opfner\\
Institut f\"ur Mathematik, Universit\"at Mainz, D--55099 Mainz\\ 
{\tt hoepfner@mathematik.uni-mainz.de}\\
{\tt http://www.mathematik.uni-mainz.de/$\sim$hoepfner}

\end{document}